\documentclass[11pt, onecolumn]{IEEEtran}
    
\usepackage{epsf, psfrag, amssymb, amsfonts, amsmath, cite,enumerate}
\usepackage{graphicx, subfigure, color,bbm}
\usepackage{boxedminipage}
\usepackage{multirow}
\usepackage{booktabs}
\usepackage{algorithm, algorithmic}
\floatstyle{plain}
\newfloat{myalgo}{tbhp}{mya}
\usepackage{setspace}
\newenvironment{Algorithm}[2][tbh]%
{\begin{myalgo}[#1]
\centering
\begin{minipage}{#2}
\begin{algorithm}[H]}%
{\end{algorithm}
\end{minipage}
\end{myalgo}}
\newtheorem{remark}{Remark}

\usepackage[scriptsize,bf]{caption}
\usepackage{etoolbox}% http://ctan.org/pkg/etoolbox
\AtBeginEnvironment{align}{\setcounter{subeqn}{0}}% Reset subequation number at start of align
\newcounter{subeqn} \renewcommand{\thesubeqn}{\theequation\alph{subeqn}}%
\newcommand{\subeqn}{%
  \refstepcounter{subeqn}% Step subequation number
  \tag{\thesubeqn}% Label equation
}
\newtheorem{theorem}{Theorem}
\newtheorem{example}{Example}

\newtheorem{lemma}{Lemma}

\newtheorem{assumption}{Assumption}
\def\psfancypar#1#2{\begingroup\def\par{\endgraf\endgroup\lineskiplimit=0pt}
               \setbox2=\hbox{\large\sc #2}
               \newdimen\tmpht \tmpht \ht2 \advance\tmpht by \baselineskip
               \font\hhuge=Times-Bold at \tmpht
               \setbox1=\hbox{{\hhuge #1}}
               \count7=\tmpht \count8=\ht1
               \divide\count8 by 1000 \divide\count7 by \count8
               \tmpht=.001\tmpht\multiply\tmpht by \count7
               \font\hhuge=Times-Bold at \tmpht
               \setbox1=\hbox{{\hhuge #1}}
               \noindent
                \hangindent1.05\wd1
               \hangafter=-2 {\hskip-\hangindent
               \lower1\ht1\hbox{\raise1.0\ht2\copy1}%
                \kern-0\wd1}\copy2\lineskiplimit=-1000pt}

\newcommand{\beq}{\begin{equation}}
\newcommand{\eeq}{\end{equation}}
\newcommand{\bqa}{\begin{eqnarray}}
\newcommand{\eqa}{\end{eqnarray}}
\newcommand{\bqn}{\begin{eqnarray*}}
\newcommand{\eqn}{\end{eqnarray*}}
\newcommand{\nn}{\nonumber}

\newcommand{\be}{\begin{enumerate}}
\newcommand{\ee}{\end{enumerate}}
\newcommand{\bi}{\begin{itemize}}
\newcommand{\ei}{\end{itemize}}
\newcommand{\bd}{\begin{description}}
\newcommand{\ed}{\end{description}}
\newcommand{\ba}{\begin{array}}
\newcommand{\ea}{\end{array}}
\newcommand{\bde}{\begin{definition}}
\newcommand{\ede}{\end{definition}}
\newcommand{\bex}{\begin{example}}
\newcommand{\eex}{\end{example}}

\def\boxit#1{\vbox{\hrule\hbox{\vrule\kern3pt
        \vbox{\kern3pt#1\kern3pt}\kern3pt\vrule}\hrule}}

\def\reals{ { {\rm  I \kern-0.15em R }  } }
\def\complex{ {\,{{\rm C} \kern-0.50em \raise0.20ex {  |}}\, }}

\def\0bf{{\bf 0}}
\def\1bf{{\bf 1}}
\def\2bf{{\bf 2}}
\def\3bf{{\bf 3}}
\def\4bf{{\bf 4}}
\def\5bf{{\bf 5}}
\def\6bf{{\bf 6}}
\def\7bf{{\bf 7}}
\def\8bf{{\bf 8}}
\def\9bf{{\bf 9}}

\def\Rbf{{\bf R}}

\def\Rxx{\Rbf_{\ssstyle X\kern-.1em X}}

\let\ssstyle=\scriptscriptstyle

\def\Kout{\setbox1=\hbox{\Huge\bf K}\hbox to
1.05\wd1{\hspace{.05\wd1}
}}
\def\Sout{\setbox1=\hbox{\Huge\bf S}\hbox to 1.05\wd1{\hspace{.05\wd1}
}}

\DeclareMathOperator*{\argmin}{arg\,min}

\begin{document}
\title{Quantized Consensus ADMM for Multi-Agent Distributed Optimization}
\author{Shengyu Zhu, Mingyi Hong, and Biao Chen}%\thanks{S. Zhu and B. Chen are with the Department of Electrical Engineering and Computer Science, Syracuse University, Syracuse, NY 13244 USA (e-mail: szhu05@syr.edu).}\thanks{M. Hong is with the }}

\maketitle

\begin{abstract}
Multi-agent distributed optimization over a network minimizes a global objective formed by a sum of local convex functions using only local computation and communication. 
We develop and analyze a quantized distributed algorithm based on the alternating direction method of multipliers (ADMM) when inter-agent communications are subject to finite capacity and other practical constraints. While existing quantized ADMM approaches only work for quadratic local objectives, the proposed algorithm can deal with more general objective functions (possibly non-smooth) including the LASSO. Under certain convexity assumptions, our algorithm converges to a consensus within $\log_{1+\eta}\Omega$ iterations, where $\eta>0$ depends on the local objectives and the network topology, and $\Omega$ is a polynomial determined by the quantization resolution, the distance between initial and optimal variable values, the local objective functions and the network topology. A tight upper bound on the consensus error is also obtained which does not depend on the size of the network.
\end{abstract}
\vspace{0.5in}

\section{Introduction}
\IEEEPARstart{T}{here} has been much research interest in distributed optimization due to recent advances in networked multi-agent systems \cite{Bertsekas:PDComputation,GG:DLWCN}. For example, {\it ad hoc} network applications may require agents to reach a consensus on the average of their measurements \cite{Ren:consensus}, including distributed coordination of mobile autonomous agents \cite{Lynch:distalgorithms}, and distributed data fusion in sensor networks \cite{Xiao:SFconsensus}. Another example is the large scale machine learning where a computation task may be executed by collaborative microprocessors with individual memories and storage spaces \cite{Bekkerman:ml1, Andrews:ml2}. Many of the distributed optimization problems, such as those mentioned, can be cast as an optimization problem of the following form
\begin{align}
\label{eqn:DistrOpt}
\min_{\tilde{x}} \sum_{i=1}^N f_i(\tilde{x}),
\end{align}
where $f_i:\mathbb{R}^M\to\mathbb{R}\cup\{\infty\}$ is the local objective function associated with agent $i$. The function $f_i$ is composed of a smooth component $g_i:\mathbb{R}^M\to\mathbb{R}\cup\{\infty\}$ and a non-smooth component $h_i:\mathbb{R}^M\to\mathbb{R}\cup\{\infty\}$, i.e., $f_i=g_i+h_i$. Examples of such models include least squares \cite{Zhu:innetworkcd, Zhu:QCviaADMM} and regularized least squares \cite{Mateos:dissparselr,Ling:Dsr,Bazerque:dss}. The variable $\tilde{x}$ may represent average temperature of a room \cite{Xiao:SFconsensus}, frequency-domain occupancy of spectra \cite{Bazerque:dss}, states of smart grid systems \cite{Gan:dprotocolforevc}, etc.  

In the above scenarios, it is commonly assumed that each agent only has the knowledge of its local objective function, and a fusion center is either disallowed or not economical. As such, the agents seek to solve (\ref{eqn:DistrOpt}) collaboratively using only local computation and communication. In practice, a number of factors, such as limited bandwidth, sensor battery power, and computing resources, place tight constraints on the rate and form of information exchange amongst neighboring agents, resulting in the {\emph{quantized communication}} constraint. The challenge of this paper is to obtain, for each agent, a reasonable solution to (\ref{eqn:DistrOpt}) in a distributed manner under the quantization constraint.

%\subsection{Related Work}
%Existing distributed approaches include , and the alternating direction method of multipliers (ADMM) \cite{Boyd:ADMM,Mota:DADMM,Wei:ADMMinDCO}. The incremental algorithm requires one to predefine a tree or loop structure in the network, whereas the subgradient descent, dual averaging, and the ADMM do not rely on any predefined structures. Subgradient descent and dual averaging are appealing due to their simplicity and the ability to handle a wide range of problems; their convergences, however, are usually slow. The decentralized ADMM approach using synchronous steps by all the agents has much faster convergence, as demonstrated by both applications \cite{Mateos:dissparselr, Ling:Dsr} and theoretic proof \cite{Wei:ADMMinDCO}. 

%While the above algorithms have been investigated extensively in the literature, we consider a practical constraint, {\em quantization}, on inter-agent communications. This is largely motivated by a number of physical factors, such as limited bandwidth, sensor battery power, and computing resources, which place tight constraints on the rate and form of information exchange amongst neighboring agents. 
Existing methods that handle this constraint include quantized incremental algorithm \cite{Rabbat:QIncremental}, quantized dual averaging \cite{Yuan:DdualQ}, and quantized subgradient method \cite{Nedic:Dsubgradient2}. The quantized incremental algorithm achieves a worst case error which is roughly $O({1}/\sqrt{k})$, where $k$ is the number of iterations, for specific quantization resolutions. This algorithm, however, does not guarantee to reach a consensus nor to converge as $k\to\infty$. The quantized dual averaging method reaches a neighborhood of the optimal solution at a rate of $O(1/\sqrt{k})$, but also does not ensure the convergence or a consensus. The quantized subgradient method converges to a consensus within a neighborhood of the optimal solution at a rate of $O(1/k)$, whose consensus error, i.e., the difference between the convergent value and the optimal value, increases in the quantization resolution, the size of the network, and the largest norm of the subgradients of local objective functions. To the best of our knowledge, there are no accelerated rates established for these algorithms when local objective functions are further known to be strongly convex. 

It is important to note that existing quantized algorithms all have sublinear rates to reach a neighborhood of the optimal solution. This is because their respective standard versions, i.e, the incremental algorithm \cite{Nedic:incrementalsub}, the dual averaging method \cite{Duchi:dualaveraging}, and the distributed subgradient descent algorithm \cite{Nedic:Dsbugradient}, have slow convergences. In addition, the errors of these quantized algorithms from the optimum tend to increase when the network becomes larger, which is much undesired as large scale networks are very typical in today's applications. The alternating direction method of multipliers (ADMM) has been known as an efficient algorithm for large scale optimizations and used in various applications such as regression and classification \cite{Boyd:ADMM}. It has been shown to have a sublinear convergence rate ${O}({1/k})$ for general convex optimization problems \cite{He:1nconvergencerate}, and to be linearly convergent for certain objective functions \cite{Hong:linearconvergence_ADMM,Deng:globallocallADMM}. Recent work of \cite{Wei:ADMMinDCO} also extends the linear convergence to a distributed ADMM method using synchronous steps when the local objective functions are strongly convex and have Lipschitz gradients. We hence expect an ADMM based quantized algorithm working well for solving (\ref{eqn:DistrOpt}) in terms of both the consensus error and the convergence time.

Unfortunately, when the quantization constraint is imposed, existing ADMM methods can only deal with quadratic local objective functions \cite{Zhu:innetworkcd, Zhu:QCviaADMM}. With dithered quantization \cite{Schuchman:ditherQ}, using the facts that quadratic functions have linear gradients and that the expectation of the dithered quantizer output is equal to the input, one can show that each agent variable converges to the optimal solution in the mean sense. For deterministic quantization, the idea is to rewrite the update as the sum a standard ADMM update of the agent variables plus an accumulated error term caused by quantization, and then use the linear convergence rate to establish convergence. These approaches, however, do not apply to general convex objective functions (see also \cite[Remark 6]{Zhu:QCviaADMM}). The local objectives can be non-smooth (e.g., the LASSO). Even when they are smooth, their gradients are not necessarily linear. Therefore, the effect of dithered quantization is hard to characterize and one can hardly write out the quantized update as the sum of a standard ADMM update plus an accumulated error term. Moreover, the linear convergence rate of the standard ADMM might fail to hold, making it more difficult to deal with quantization.

Our main contribution is to develop a quantized distributed ADMM algorithm using deterministic quantization. We do not directly prove the convergence of the variables at each agent; instead, we seek to establish the convergence of an auxiliary vector which determines the update of the agent variables. In particular, we show that this algorithm converges to a consensus within finite iterations under certain convexity assumptions as long as an initialization condition is satisfied. The initialization condition is rather mild; indeed, simply setting all the variables to $0$ suffices. We derive a tight upper bound on the consensus error which does not depend on the size of the network. We finally characterize the convergence time, that is, our algorithm converges within $\log_{1+\eta}\Omega$ iterations where $\eta>0$ depends on the local objectives and the network topology, and $\Omega$ is a polynomial decided by the quantization resolution, the initial variable values, the local objective functions and the network topology. The proof idea also provides a framework for convergence proof of other quantized algorithms (see Section \ref{sec:convergence}).

%By comparison to our previous work \cite{Zhu:QCviaADMM} where the local objective functions are quadratic, the problem setting and proof of this paper are more general. Indeed, the general local objective functions make the problem much harder because of the possible nonlinearity of their gradients: the effect of dithered quantization (see, e.g., \cite{Schuchman:ditherQ}) is not easy to characterize due to this nonlinearity. We thereby adopt deterministic quantization in this paper. Another concern is that one can hardly write out the quantized update as the sum of an ideal distributed ADMM update plus an accumulated error term. As such, we do not directly prove the convergence of the variables at each node; instead, we seek to show that an auxiliary vector which determines the updates of variables at each agent converges. Our idea also provides a framework for convergence proof of other quantized algorithms (see Section \ref{sec:convergence}).% and Remark \ref{eqn:PfIdeaGeneralization}).

The rest of this paper is organized as follows. Section \ref{sec:DOviaADMM} reviews the application of the ADMM to distributed optimization without the quantization constraint, resulting in a distributed ADMM algorithm. In Section \ref{sec:QCADMM}, we use deterministic quantization to modify the above algorithm to handle the quantization constraint. We show the relation of the quantized algorithm to the standard ADMM and establish the desired convergence results. Simulations are provided in Section \ref{sec:simulations}, followed by conclusion in Section \ref{sec:conclusion} along with discussions on future research directions.

{\it Notations:} 
We use $0$ to denote the all-zero column vector with a suitably defined dimension. $1_K$ is the $K$-dimensional all-one column vector; $0_K$ and $I_K$ are the $K\times K$ all-zero and identity matrix, respectively. Notation $\otimes$ denotes the Kronecker product and $\|x\|_2$ denotes the Euclidean norm of a vector $x$. Given a positive semidefinite matrix $G$ with proper dimensions, the $G$-norm of $x$ is $\|x\|_G=\sqrt{x^TGx}$. Denote $\sigma_{\max}(D)$ as the largest singular value of a matrix $D$ and $\tilde{\sigma}_{\min}(D)$ as the smallest nonzero singular value of $D$. $\partial f(x)$ denotes a subgradient of $f$ at $x$ for a convex function $f(x)$ while $\nabla f(x)$ denotes the gradient if it is known to be differentiable.

We use two definitions of rate of convergence for an iterative algorithm. A sequence $x^k$, where the superscript $k$ stands for time index, is said to converge \emph{Q-linearly} to a point $x^*$ if there exists a number $\upsilon\in(0,1)$ such that $\lim_{k\to\infty}\frac{\|x^{k+1}-x^*\|}{\|x^k-x^*\|}=\upsilon$ with $\|\cdot\|$ being a vector norm. We say that a sequence $y^k$ converges \emph{R-linearly} to $y^*$ if $\|y^k-y^*\|\leq\|x^k-x^*\|$ for all $k$, where $x^k$ converges Q-linearly to $x^*$.
%\subsection{Related Work}
%\subsection{Our Contributions}
%\subsection{Paper Organization}
%\subsection{Notations}
\section{Distributed Optimization via the ADMM}
\label{sec:DOviaADMM}
This section reviews the consensus ADMM (C-ADMM) for distributed optimization where agents can send and receive real data with infinite precision. This ideal case provides a good understanding of how the ADMM works and performs in a distributed manner. We start with the problem setting and assumptions.
\label{sec:consensusADMM}
\subsection{Problem Setting and Assumptions}
Throughout the paper we consider a network consisting of $N$ agents bidirectionally connected by $E$ edges, where each agent $i$ has its own objective function $f_i:\mathbb{R}^M\to\mathbb{R}\cup\{\infty\}$. Assume that the network topology is fixed. We describe this network as a symmetric directed graph $\mathcal{G}_d=\{\mathcal{V},\mathcal{A}\}$ or an undirected graph $\mathcal{G}_u=\{\mathcal{V},\mathcal{E}\}$, where $\mathcal{V}$ is the set of vertices with cardinality $|\mathcal{V}|=N$, $\mathcal{A}$ is the set of arcs with $|\mathcal{A}|=2E$, and $\mathcal{E}$ is the set of edges with $|\mathcal{E}|=E$. Based on this graph, we would like to develop in-network algorithms that find the global optimum $\tilde{x}^*$ (not necessarily unique) minimizing $$\sum_{i=1}^N f_i(\tilde{x}).$$

We make the following assumptions on the local objective functions $f_i=g_i+h_i,i=1,2,\cdots,N$.
\begin{assumption}
\label{asm:convex}
The local objective functions are proper closed convex functions; for every $\tilde{x}$ where $f_i(\tilde{x})$ is well defined and $f_i(\tilde{x})<\infty$, there exists at least one bounded subgradient $\partial f_i(\tilde{x})$ such that $$f_i(\tilde{y})\geq f_i(\tilde{x})+\left(\partial f(\tilde{x})\right)^T(\tilde{y}-\tilde{x}), \forall\tilde{y}\in\mathbb{R}^M.$$ Moreover, the minimum of (\ref{eqn:DistrOpt}) can be attained. 
\end{assumption}
\begin{assumption}
\label{asm:StrLipconvex}
The smooth components have Lipschitz continuous gradients, i.e., for each agent $i$ there exists some $M_{g_i}>0$ such that $$\|\nabla g_i(\tilde{x})-\nabla g_i(\tilde{y})\|_2 \leq M_{g_i}\|\tilde{x}-\tilde{y}\|_2, \forall \tilde{x},\tilde{y}\in\mathbb{R}^M.$$
In addition, the smooth components are strongly convex, i.e., for each agent $i$ there exists some $m_{g_i}>0$ such that $$\left(\nabla g_i(\tilde{x})-\nabla g_i(\tilde{y})\right)^T(\tilde{x}-\tilde{y})\geq m_{g_i}\|\tilde{x}-\tilde{y}\|_2^2, \forall \tilde{x},\tilde{y}\in\mathbb{R}^M.$$
\end{assumption}
\begin{assumption}
\label{asm:nonsmooth}
The non-smooth components $h_i$'s are convex.
\end{assumption}
%\begin{assumption}
%\label{asm:lineargradient}
%The local objective functions have linear gradients, i.e., for each agent $i$ and $u,v\in\mathbb{R}^M$, $\nabla f_i(u+v)=\nabla f_i(u)+\nabla f_i(v)$.
%\end{assumption}

Note that Assumption \ref{asm:StrLipconvex} implies the differentiability of $g_i$. Assumptions \ref{asm:convex}--\ref{asm:nonsmooth} together indicate that (\ref{eqn:DistrOpt}) has a unique and attainable solution, i.e., $\tilde{x}^*\in\mathbb{R}^M$ is unique. We only need Assumption \ref{asm:convex} to show the convergence of the C-ADMM, while Assumption \ref{asm:StrLipconvex} is essential to establish the linear convergence when $f_i$ only contains the smooth component $g_i$.%, which is the key to proving our main result in Section \ref{sec:QCADMM}.

\subsection{The ADMM for Distributed Optimization: C-ADMM}
To solve (\ref{eqn:DistrOpt}) using the ADMM, we first reformulate it as 
\begin{equation}
\begin{aligned}
\label{eqn:admmformulation}
& \underset{\{x_i\},\{z_{ij}\}}{\text{minimize}}
& & \sum_{i=1}^N f_i(x_i)\\
& \text{subject to}
& & x_i=z_{ij},x_j=z_{ij},\forall (i,j)\in\mathcal{A},
\end{aligned}
\end{equation}
where $x_i\in\mathbb{R}^M$ is the local copy of the common optimization variable $\tilde{x}$ at agent $i$ and $z_{ij}\in\mathbb{R}^M$ is an auxiliary variable imposing the consensus constraint on neighboring agents $i$ and $j$. As the given network is connected, the consensus constraint ensures the consensus to be achieved over the entire network, i.e., $x_i=x_j, \forall i,j \in \mathcal{A}$, which in turn guarantees that (\ref{eqn:admmformulation}) is equivalent to (\ref{eqn:DistrOpt}). Further define $x\in\mathbb{R}^{NM}$ as a vector concatenating all $x_i$, $z\in\mathbb{R}^{2EM}$ as a vector concatenating all $z_{ij}$, $g(x)=\sum_{i=1}^N g_i(x_i)$,  $h(x)=\sum_{i=1}^N h_i(x_i)$, and $f(x)=g(x)+h(x)$. Then (\ref{eqn:admmformulation}) can be written in a matrix form as 
\begin{equation}
\label{eqn:matrixform}
\begin{aligned}
& \underset{x,z}{\text{minimize}}
& & f(x)+{f}'(z)\\
& \text{subject to}
& & Ax+Bz=0,
\end{aligned}
\end{equation}
where ${f}'(z)=0$, $B = [-I_{2EM}; -I_{2EM}]$ with $I_{2EM}$ being a $2EM\times 2EM$ identity matrix, and $A = [A_1;A_2]$ with $A_1,A_2\in\mathbb{R}^{2EM\times NM}$ both being composed of $2E\times M$ blocks of $M\times M$ matrices. If $(i,j)\in\mathcal{A}$ and $z_{ij}$ is the the $q$th block of $z$, then the $(q,i)$th block of $A_1$ and the $(q,j)$th block of $A_2$ are $M\times M$ identity matrices $I_M$; otherwise the corresponding blocks are $M\times M$ zero matrices $0_M$. 

We are now ready to apply the ADMM to solving (\ref{eqn:DistrOpt}). The augmented Lagrangian of (\ref{eqn:matrixform}) is
\begin{align}
\label{eqn:Lagragian}
L_\rho(x,z,\lambda) =f(x)+\lambda^T(Ax+Bz)+\frac{\rho}{2}\|Ax+Bz\|_2^2,
\end{align}
where $\lambda=[\beta;\gamma]$ with $\beta,\gamma\in\mathbb{R}^{2EM}$ is the Lagrange multiplier and $\rho\in\mathbb{R}$ is a positive algorithm parameter. At iteration $k+1$, the ADMM first obtains $x^{k+1}$ by minimizing $L_\rho(x,z^k,\lambda^k)$, then calculates $z^{k+1}$ by minimizing $L_\rho(x^{k+1},z,\lambda^k)$, and finally updates $\lambda^{k+1}$ from $x^{k+1}$ and $z^{k+1}$. Noting that both $A$ and $B$ are full column-rank for connected networks, Assumption \ref{asm:convex} implies that such $x^{k+1}$ and $z^{k+1}$ exist uniquely. We then have the ADMM update given by
\begin{equation}
\label{eqn:admmupdates}
\begin{aligned}
x\text{-update:}&~\partial f(x^{k+1}) + A^T\lambda^k+\rho A^T(Ax^{k+1}+Bz^k)=0,\\
z\text{-update:}&~B^T\lambda^k+\rho B^T(Ax^{k+1}+Bz^{k+1}) = 0,\\
\lambda\text{-update:}&~\lambda^{k+1}-\lambda^k-\rho(Ax^{k+1}+Bz^{k+1})=0,
\end{aligned}
\end{equation} 
where $\partial f(x^{k+1})$ denotes a subgradient of $f(x)$ at $x = x^{k+1}$.

A nice convergence property of the ADMM, known as \emph{global convergence}, states that the sequence $( x^k, z^k,\lambda^k)$ generated by (\ref{eqn:admmupdates}) has a single limit point  $(x^*,z^*,\lambda^*)$ under Assumption \ref{asm:convex}. Proofs can be found in \cite{Boyd:ADMM,Deng:globallocallADMM,He:1nconvergencerate}. If $\tilde{x}^*$ is unique (e.g., when Assumptions \ref{asm:convex}--\ref{asm:nonsmooth} hold), then $x^*={\bf 1}_N\tilde{x}^*$ is also unique where ${\bf 1}_N=1_N\otimes I_M$, i.e., a matrix consisting of $N\times1$ blocks of $I_M$. To summarize, we have
\begin{lemma}[Global convergence of the ADMM \cite{Boyd:ADMM,He:1nconvergencerate,Deng:globallocallADMM}]
\label{lem:globalconvergence}
Under Assumption \ref{asm:convex}, the updates in (\ref{eqn:admmupdates}) yield that for any initial values $x^0\in\mathbb{R}^{NM}$, $z^0\in\mathbb{R}^{2EM}$ and $ \lambda^0\in\mathbb{R}^{4EM}$, $$x^k\to x^*,~z^k\to z^*,~\text{and}~ \lambda^k\to\lambda^*~\text{as}~k\to\infty,$$
where $(x^*,z^*,\lambda^*)$ is a primal-dual solution to (\ref{eqn:Lagragian}). If Assumptions \ref{asm:StrLipconvex} and \ref{asm:nonsmooth} also hold, then $x^*={\bf 1}_N\tilde{x}^*$ is the unique solution to (\ref{eqn:matrixform}).
\end{lemma}

While (\ref{eqn:admmupdates}) provides an efficient {centralized} algorithm to solve (\ref{eqn:matrixform}), it is not clear whether it can be carried out in a distributed manner, i.e., data exchanges only occur among neighboring nodes. Interestingly, as established in Lemma \ref{lem:globalconvergence}, the ADMM allows any initial values $x^0,z^0$ and $ \lambda^0$ for global convergence under Assumption \ref{asm:convex}; there indeed exist initial values that decentralize (\ref{eqn:admmupdates}). Define $M_+=A_1^T+A_2^T$ and $ M_-=A_1^T- A_2^T$, which are respectively the extended unoriented and oriented incident matrices with respect to the directed graph $\mathcal{G}_d$.\footnote{By ``extended", we mean the Kronecker product of the original matrix multiplying $I_M$; see \cite{Chung:spercalG, Fiedler:connectivityG, Cvetkovic:LaplaciansFG}.} As shown in \cite{Wei:ADMMinDCO}, by initializing $\beta^0= -\gamma^0$ and $z^0=\frac{1}{2} M_+^T x^0$, the update in (\ref{eqn:admmupdates}) leads to 
\begin{equation}
\label{eqn:distributedversion}
\begin{aligned}
x_i\text{-update:}&~\partial f(x_i^{k+1}) + 2\rho|\mathcal{N}_i|x_i^{k+1}+\alpha_i^k -\rho\left(|\mathcal{N}_i|x_i^k+\sum_{j\in\mathcal{N}_i} x_j^k\right)=0,\\
%&x_i^{k+1}=(\partial f_i+2\rho|\mathcal{N}_i|I_M)^{-1}\left(\rho|\mathcal{N}_i|x_i^k+\rho\sum_{j\in\mathcal{N}_i} x_j^k-\alpha_i^k\right),\\
\alpha_i\text{-update:}&~\alpha_i^{k+1}=\alpha_i^k+\rho\left(|\mathcal{N}_i|x_i^{k+1}-\sum_{j\in\mathcal{N}_i} x_j^{k+1}\right)
\end{aligned}
\end{equation} 
at node $i$, where $\mathcal{N}_i$ denotes the set of neighbors of node $i$ and $\alpha_i^k\in\mathbb{R}^M$ is the $i$th block of $\alpha^k=M_-\beta^k\in\mathbb{R}^{NM}$. Obviously, (\ref{eqn:distributedversion}) is fully decentralized as the update of $x_i^{k+1}$ and $\alpha_i^{k+1}$ only relies on local and neighboring information. We refer to (\ref{eqn:distributedversion}) as the C-ADMM update.
\subsection{Linear Convergence of the C-ADMM}
\label{sec:lcCADMM}
Before stating the convergence results of the C-ADMM, we introduce some useful facts that are related to the undirected graph $\mathcal{G}_u$. These facts not only simplify our presentation but also help establish the main theorem in Section \ref{sec:convergence} .

Define $L_+=\frac{1}{2}M_+M_+^T$ and $L_-=\frac{1}{2}M_-M_-^T$, which are respectively the extended signless and signed Laplacian matrices with respect to $\mathcal{G}_u$. Then $W=\frac{1}{2}(L_++L_-)$ is the extended degree matrix, i.e., a block diagonal matrix with its $(i,i)$th block being the Kronecker product of $|\mathcal{N}_i|$ multiplying $I_M$ and other blocks being $0_M$. We have the following lemma regarding $L_-$.
\begin{lemma}[\cite{Chung:spercalG, Fiedler:connectivityG, Cvetkovic:LaplaciansFG}]
\label{lem:LminusP}
For connected networks, $L_-$ is positive semidefinite and always has $0$ as its eigenvalue; $L_-b=0$ if and only if $b={\bf 1}_N\tilde{b}$ for some $\tilde{b}\in\mathbb{R}^M$.
\end{lemma}

	As a result of the above lemma, we obtain Lemma \ref{lem:abrelation} which states the one-to-one correspondence between $\alpha$ and $\beta$ provided that $\beta$ lies in the column space of $M_-^T$. 
\begin{lemma}
	\label{lem:abrelation}
	Given a connected network, if $\beta$ lies in the column space of $M_-^T$, then $\alpha$ and $\beta$ are one-to-one correspondence; i.e., let $\alpha=M_-\beta$ and $\alpha'=M_-\beta'$ for some $\beta$ and $\beta'$ in the column space of $M_-^T$, then $\alpha=\alpha'$ if and only if $\beta=\beta'$. 
\end{lemma}
\begin{IEEEproof}
	That $\beta=\beta'$ implies $\alpha=\alpha'$ is straightforward. Consider $\alpha=M_-\beta$ and write $\beta= M_-^Tb$  for some $b\in\mathbb{R}^{NM}$. $\alpha'$, $\beta'$ and $b'$ are similarly defined. Then
	$$\alpha-\alpha' = 2L_-(b-b')=0,$$
which implies $b-b'={\bf 1}_N\tilde{b}$ for some $\tilde{b}\in\mathbb{R}^M$ from Lemma \ref{lem:LminusP}. Since $M_-$ is the extended oriented incident matrix with respect to $\mathcal{G}_d$, we have
$$\beta-\beta' = M_-^T(b-b') =\left(M_-^T{\bf 1}_N\right)\tilde{b} = 0.$$
\end{IEEEproof}
%With the notions of $L_-$, $L_+$ and $W$, (\ref{eqn:distributedversion}) can be written in a matrix form as
%\begin{align}
%\label{eqn:compactform2}
 %\left[\begin{array}{c} x^{k+1} \\\alpha^{k+1}\end{array} \right] =\left[\begin{matrix} \nabla f + 2\rho W & 0\\ -\rho L_- & I \end{matrix}\right ] ^{-1} \left[\begin{matrix} \rho L_+ & -I \\ 0 & I\end{matrix}\right] \left[\begin{array}{c} x^{k} \\\alpha^{k}\end{array} \right],
%\end{align} 
%where $\nabla f + 2W$ shall be understood as a $NM\times NM$ matrix function that maps $\mathbb{R}^{NM}\to\mathbb{R}^{NM}$, and its inverse exists due to Assumption \ref{asm:convex}. 

We now turn our attention to the convergence properties of the C-ADMM. That the C-ADMM converges follows directly from global convergence of the ADMM [cf. Lemma \ref{lem:globalconvergence}]. To establish its linear convergence, we state the following lemma regarding the convergence rate of a vector concatenating $z$ and $\beta$.\footnote{The linear convergence results of \cite{Hong:linearconvergence_ADMM,Deng:globallocallADMM} do not apply here. The step size of the dual variable update need be sufficiently small in \cite{Hong:linearconvergence_ADMM} while the C-ADMM has a fixed step size $\rho$. The linear convergence result in \cite{Deng:globallocallADMM} requires that $f'(z)$ is strongly convex or $B$ is full row-rank. However, in our formulation (\ref{eqn:matrixform}), $f'(z)=0$ is not strongly convex and $B = [-I_{2EM}; -I_{2EM}]$ is row-rank deficient.}

\begin{lemma}[{\cite[Theorem 1]{Wei:ADMMinDCO}}]
	\label{lem:linearconvergence} 
	Consider the ADMM iteration (\ref{eqn:admmupdates}) that solves (\ref{eqn:matrixform}). Define $$u = \left[\begin{array}{c} z\\ \beta \end{array} \right]~\text{and}~G = \begin{bmatrix} \rho I_{2EM} & 0_{2EM} \\ 0_{2EM}& \frac{1}{\rho} I_{2EM}  \end{bmatrix} $$ with $\beta$ being a dual variable. Let $z^0=\frac{1}{2}M_+^Tx^0$ and $\beta^0=-\gamma^0$, where $\gamma$ is the other dual variable and $\beta^0$ is initialized to lie in the column space of $M_-^T$. Suppose that Assumptions \ref{asm:convex}--\ref{asm:nonsmooth} hold. Then
	%\begin{enumerate}
	%\item under Assumption \ref{asm:convex}, $z^k = \frac{1}{2}M_+^Tx^k$ and $\beta^k=-\gamma^k$ lying in the column space of $M_-^T$ for $k =0,1,\cdots$
	$z^k = \frac{1}{2}M_+^Tx^k$ and $\beta^k=-\gamma^k$ lying in the column space of $M_-^T$ for $k =0,1,\cdots$, and $(x^k,z^k,\beta^k)$ converges uniquely to $(x^*,z^*,\beta^*)$ where $x^*={\bf 1}_N\tilde{x}^*$, $z^*=\frac{1}{2}M_+^Tx^*$ and $\beta^*$ is a vector in the column space of $M_-^T$. If we further have $f_i = g_i$, then for any $\mu > 1$, $u^k = [z^k;\beta^k]$ converges Q-linearly to its optimal $u^*=[z^*;\beta^*]$ with respect to the $G$-norm 
	\begin{equation}
	\label{eqn:Ulinearconvergence}
	\begin{aligned}
	\|u^{k+1}-u^*\|_G\leq\frac{1}{1+\eta}\|u^k-u^*\|_G,
	\end{aligned}
	\end{equation}
	where $\eta=\sqrt{1+\delta}-1$, $m_g\triangleq\min_i\{m_{g_i}\}$, $M_g\triangleq\max_i\{M_{g_i}\}$, and 
\begin{align}
\delta &= \min \left\{\frac{(\mu-1)\tilde{\sigma}^2_{\min}(M_-)}{\mu \sigma^2_{\max}(M_+)},\frac{4\rho m_g\tilde{\sigma}_{\min}^2(M_-)}{\rho^2\sigma^2_{\max}(M_+)\tilde{\sigma}^2_{\min}(M_-)+\mu M_g^2}\right\}.\nn
\end{align}
%Let $m_f = \min_i\{m_{f_i}\}$. Then $x^k$ converges R-linearly to $x^*$ due to
%\begin{align}
%\label{eqn:xk1uk}
%\|x^{k+1}-x^*\|_2^2\leq\frac{1}{m_f}\|u^k-u^*\|_G^2,\nn
%\end{align}
%\end{enumerate}
\end{lemma}
See \cite{Wei:ADMMinDCO} for the proof.

Notice that when $f_i=g_i$, the non-smooth component $h_i=0$ and Assumption \ref{asm:nonsmooth} is satisfied automatically. Then the above theorem indicates that $u^k$ is linearly convergent to the unique optimum provided that the ADMM update is initialized properly, $f_i$ is smooth, and Assumptions  \ref{asm:convex} and \ref{asm:StrLipconvex} hold. As an interesting observation from the ADMM iteration (\ref{eqn:admmupdates}), we notice that $(x^{k+1},z^{k+1},\lambda^{k+1})$ is obtained only based on $(z^k, \lambda^k)$. Given the initialization in Lemma \ref{lem:linearconvergence}, $\lambda^k=[\beta^k;-\beta^k]$ for $k=0,1,\cdots$, and hence the $(k+1)$th update only requires the knowledge of $u^k=[z^k;\beta^k]$. We can then define a function $\psi(\cdot)$, which represents the update of $u^{k+1}$ from $u^k$ via (\ref{eqn:admmupdates}), as 
\begin{align}
\label{eqn:uh}
u^{k+1}=\psi(u^k).
\end{align}
Then (\ref{eqn:Ulinearconvergence}) is equivalent to
\begin{align}
\label{eqn:uhlc}
\|\psi(u^k)-u^*\|_G\leq\frac{1}{1+\eta}\|u^k-u^*\|_G.
\end{align}

Compared with the initialization conditions that lead the ADMM to the C-ADMM, we notice an extra initialization condition that $\beta_0$ lies in the column space of $M_-^T$ in Lemma \ref{lem:linearconvergence}. Without this condition, $x^k$ still converges to its unique optimal value $x^*={\bf{1}}_N\tilde{x}^*$ under Assumptions \ref{asm:convex}--\ref{asm:nonsmooth}, whereas the uniqueness of $\beta^*$ may not hold nor the linear convergence of $u^k$. While this condition implies that $\alpha^0$ in the C-ADMM is initialized in the column space of $L_-$, it is not clear whether such $\alpha^0$ ensures the linear convergence of the C-ADMM. In the following we show that the C-ADMM indeed converges linearly to its optimal value under Assumptions  \ref{asm:convex} and \ref{asm:StrLipconvex} given that $f_i=g_i$ and $\alpha^0$ lies in the column space of $L_-$. To see this, let $x^0$ and $\alpha^0$ be the initial values in the C-ADMM where $\alpha^0$ is in the column space of $L_-$. Lemma \ref{lem:abrelation} implies that there exists a unique $\beta^0$ in the column space of $M_-^T$ such that $\alpha^0=M_-\beta^0$. This $\beta^0$ together with $z^0=\frac{1}{2}M_+^Tx^0$ and $\gamma^0=-\beta^0$ leads the ADMM to the C-ADMM with exactly the same initial values $x^0$ and $\alpha^0$. Then $x^k$ stays the same in the ADMM and the C-ADMM iterations, and $\alpha^k=M_-\beta^k$ where $\beta^k$ is the $k$-th update in the ADMM. As such, we can study the convergence rate of the C-ADMM through the ADMM that is correspondingly initialized. The linear convergence of the C-ADMM is presented below.
\begin{theorem}[Linear convergence of the C-ADMM]
\label{thm:lcC-ADMM}
Let $f_i=g_i$ and consider the C-ADMM update (\ref{eqn:distributedversion}). Suppose that Assumptions \ref{asm:convex} and \ref{asm:StrLipconvex} hold. If we initialize $\alpha^0$ in the column space of $L_-$, then $[x^k;\alpha^k]$ converges R-linearly to its optimal value $[x^*;\alpha^*]$ where $x^*={\bf 1}_N\tilde{x}^*$ and $\alpha^*= \nabla g({\bf 1}_N\tilde{x}^*)$.
\end{theorem}
\begin{IEEEproof}
From Lemma \ref{lem:abrelation}, we can denote by $\beta^0$ the unique vector in the column space of $M_-^T$ such that $\alpha^0=M_-\beta_0$. Consider the ADMM iteration (\ref{eqn:admmupdates}) with initialization $z^0=\frac{1}{2}M_+^Tx^0$ and $\lambda^0 = [\beta^0;-\beta^0]$. Then the convergence of $x^k$ to $x^*={\bf 1}_N\tilde{x}^*$ follows from Lemma \ref{lem:globalconvergence}. By plugging $x_i^*=\tilde{x}^*$ into (\ref{eqn:distributedversion}), we obtain $\alpha_i^*=\nabla g_i(\tilde{x}^*)$.

To show the linear convergence, we first have from \cite[Equation (29)]{Wei:ADMMinDCO} that $$\|x^{k+1}-x^*\|_2^2\leq\frac{1}{m_g}\|u^k-u^*\|_G^2,$$ where $m_g$ and $u^k$ are defined in Lemma \ref{lem:linearconvergence}. Recalling the definition of $G$, we have 
\begin{align}
\|\alpha^{k+1}-\alpha^*\|_2^2&=\|M_-(\beta^{k+1}-\beta^*)\|^2_2\nn\\
&\leq\sigma^2_{\max}(M_-)\rho\|u^{k+1}-u^*\|_G^2\nn\\
&\stackrel{(a)}{\leq}\sigma^2_{\max}(M_-){\frac{\rho}{(1+\eta)^2}}\|u^k-u^*\|_G^2,\nn
\end{align}
where $(a)$ is from the linear convergence of $u^k$ as the initialization conditions in Lemma \ref{lem:linearconvergence} are satisfied. Therefore,
$$\left\|\left[\begin{array}{c} x^{k+1} \\\alpha^{k+1}\end{array} \right] -\left[\begin{array}{c} x^* \\\alpha^*\end{array} \right]\right\|_2^2\leq \left(\frac{1}{m_g}+\sigma^2_{\max}(M_-){\frac{\rho}{(1+\eta)^2}}\right)\|u^k-u^*\|_G^2,$$
which establishes the R-linear convergence of $[x^k;\alpha^k]$ since $u^k$ converges Q-linearly to $u^*$.
\end{IEEEproof}

\section{Quantized Consensus ADMM}
\label{sec:QCADMM}
%Quantization appears in many real applications due to the physical constraints. A useful approach is to add dither to the quantizer, which makes the quantization error sequence is i.i.d. and independent of the input signal sequence; see, e.g., . The idea is to utilize the linearity of the expectation operator; however, this property does not apply here as the gradients or subgradients of local objective functions are not always linear. Therefore, we only consider the deterministic quantization. 
To model the effect of quantized communication, we assume that each agent can store and  compute real values with infinite precision; an agent, however, can only transmit quantized data through the channel which are received by its neighbors without any error. Given a quantization resolution $\Delta>0$, define the quantization lattice in $\mathbb{R}$ by $$\Lambda = \{t\Delta:t\in\mathbb{Z}\}.$$ A quantizer is a function $Q:\mathbb{R}\to\Lambda$ that maps a real value to some point in $\Lambda$. Among all deterministic quantizers, we consider the rounding quantizer that projects $y\in\mathbb{R}$ to its nearest point in $\Lambda$: 
\begin{equation}
\label{eqn:Qr}
\begin{aligned}
Q(y) = t \Delta,~\text{if}~\left(t-\frac{1}{2}\right)\Delta\leq y< \left(t+\frac{1}{2}\right)\Delta.
\end{aligned}
\end{equation}
By quantizing a vector we mean quantizing each of its entries. For $w\in\mathbb{R}^L,L\in\mathbb{Z}^+$, the rounding quantizer projects $w$ to its nearest point in $\Lambda^L$, and we use $w_{[Q]}$ to denote the quantizer output of $w$. Define $e=w_{[Q]}-w$ as the quantization error. It is clear that for any $w\in\mathbb{R}^L$,
\begin{align}
\label{eqn:qerr}
\|e\|_2\leq\frac{1}{2}{\Delta}\sqrt{L}.
\end{align}

We next use the above rounding quantization to modify the C-ADMM to meet the communication constraint, resulting in the quantized consensus ADMM (QC-ADMM) in Algorithm \ref{tab:QCADMM}.
\begin{Algorithm}{.9\textwidth}
	\caption{QC-ADMM for solving (\ref{eqn:DistrOpt})}
	%\begin{minipage}{0.8\textwidth}
	\begin{algorithmic}[1]
	\label{tab:QCADMM}
	\REQUIRE Initialize~$x_i^0\in\mathbb{R}^M$ and $\alpha_{iQ}^0\in\mathbb{R}^M$ for each agent $i,i=1,2,\cdots,N$ such that $\alpha_Q^0$ lies in the column space of $L_-$. Set $\rho>0$ and $k=0$.
		%\STATE {\bf Set} $k=1.$
	%\FOR{$k=0,1,\cdots,K-1$, every agent $i$}
	%\STATE \begin{align} x_i^{k+1}&\gets(\nabla f_i+2\rho|\mathcal{N}_i|I_M)^{-1}\left(\rho|\mathcal{N}_i|x_{i[Q]}^k+\rho\sum_{j\in\mathcal{N}_i} x_{j[Q]}^k-\alpha_{iQ}^k\right),\nn\\
%\alpha_{iQ}^{k+1}&\gets\alpha_{iQ}^k+\rho\left(|\mathcal{N}_i|x_{i[Q]}^{k+1}-\sum_{j\in\mathcal{N}_i} x_{j[Q]}^{k+1}\right).\nn
	%\end{align}
	%\ENDFOR
	%\STATE {\bf set} $x^0=x_{[Q]}^K$, $\alpha_Q^0=0$ and $k=0$.
	\REPEAT
			\STATE every agent $i$ {\bf do}
			\begin{equation} 
			\label{eqn:distributedversionQ}
			\begin{aligned}
			x_i\text{-update:}&~\partial f_i(x_i^{k+1})+2\rho|\mathcal{N}_i|x_i^{k+1}+\alpha_{iQ}^k-\rho\left(|\mathcal{N}_i|x_{i[Q]}^k+\sum_{j\in\mathcal{N}_i} x_{j[Q]}^k\right)=0,\\
\alpha_{iQ}\text{-update:}&~\alpha_{iQ}^{k+1}=\alpha_{iQ}^k+\rho\left(|\mathcal{N}_i|x_{i[Q]}^{k+1}-\sum_{j\in\mathcal{N}_i} x_{j[Q]}^{k+1}\right).
			\end{aligned}
			\end{equation}
			\STATE {\bf set} $k=k+1$.	
	\UNTIL{a predefined stopping criterion (e.g., a maximum iteration number) is satisfied.}
 %(e.g., a maximum iteration number) }%
	\end{algorithmic}
	%\end{minipage}
\end{Algorithm}
%\begin{equation}
%\label{eqn:distributedversionQ}
%\begin{aligned}
%&x_i^{k+1}=(\nabla f_i+2\rho|\mathcal{N}_i|I_M)^{-1}\left(\rho|\mathcal{N}_i|x_{i[Q]}^k+\rho\sum_{j\in\mathcal{N}_i} x_{j[Q]}^k-\alpha_{iQ}^k\right),\\
%&\alpha_{iQ}^{k+1}=\alpha_{iQ}^k+\rho\left(|\mathcal{N}_i|x_{i[Q]}^{k+1}-\sum_{j\in\mathcal{N}_i} x_{j[Q]}^{k+1}\right).
%\end{aligned}
%\end{equation}

In Algorithm \ref{tab:QCADMM} we use the subscript $Q$ to differentiate between the QC-ADMM and C-ADMM updates, and $\alpha_{iQ}^k$ is not necessarily equal to $\alpha^k_{i[Q]}$. Note that $x_i^{k}$ is quantized at its own node for the $(k+1)$th update; the reason will be given in Remark \ref{rmk:otherscheme}. We will establish the connection of the QC-ADMM with the standard ADMM in Section \ref{sec:connection}, and the convergence results of the QC-ADMM in Sections \ref{sec:convergence} and \ref{sec:convergenceG}. 

\subsection{Connection with the ADMM}
\label{sec:connection}
Now the QC-ADMM update (\ref{eqn:distributedversionQ}) seems to be a direct modification from the C-ADMM update by quantizing $x^k$ for the $(k+1)$th update; it is not clear how it relates to the standard ADMM. We will show that (\ref{eqn:distributedversionQ}) can be derived from (\ref{eqn:admmupdates}) by imposing a {{quantization}} operation on $x$ immediately after the $x$-update, i.e., 
\begin{equation}
\label{eqn:admmupdatesQ}
\begin{aligned}
x\text{-update}:&~~\partial f(x^{k+1}) + A^T\lambda_{Q}^k+\rho A^T(Ax^{k+1}+Bz_{Q}^k)=0,\\
x_{[Q]}\text{-update}:&~~ x^{k+1}_{[Q]}=Q(x^{k+1}),\\
z_{Q}\text{-update}:&~~B^T\lambda_{Q}^k+\rho B^T(Ax_{[Q]}^{k+1}+Bz_{Q}^{k+1}) = 0,\\
\lambda_{Q}\text{-update}:&~~\lambda_{Q}^{k+1}-\lambda_{Q}^k-\rho(Ax_{[Q]}^{k+1}+Bz_{Q}^{k+1})=0,
\end{aligned}
\end{equation} 
where the subscript $Q$ is adopted to differentiate between the updates before and after the $x_{[Q]}$-update. Again we do not have $z^k_{Q}=z^k_{[Q]}$ or $\lambda^k_{Q}=\lambda^k_{[Q]}$ in general. Since $(x^{k+1}, z^{k+1}_Q, \lambda^{k+1}_Q)$ is updated only based on $(z^k_{Q}, \lambda^k_{Q})$ given the deterministic quantization operation defined by (\ref{eqn:Qr}), we can still perform the $\lambda$-update and $z$-update before the $x_{[Q]}$-update. That is, (\ref{eqn:admmupdatesQ}) is equivalent to
\begin{align}
x\text{-update}:&~~\partial f(x^{k+1}) + A^T\lambda_{Q}^k+\rho A^T(Ax^{k+1}+Bz_{Q}^k)=0\refstepcounter{equation}\subeqn\label{eqn:aa},\\
z\text{-update}:&~~B^T\lambda_{Q}^k+\rho B^T(Ax^{k+1}+Bz^{k+1}) = 0\subeqn\label{eqn:ab},\\
\lambda\text{-update}:&~~\lambda^{k+1}-\lambda_{Q}^k-\rho(Ax^{k+1}+Bz^{k+1})=0\subeqn\label{eqn:ac},\\
x_{[Q]}\text{-update}:&~~ x^{k+1}_{[Q]}=Q(x^{k+1})\subeqn\label{eqn:ad},\\
z_{Q}\text{-update}:&~~B^T\lambda_{Q}^k+\rho B^T(Ax_{[Q]}^{k+1}+Bz_{Q}^{k+1}) = 0\subeqn\label{eqn:ae},\\
\lambda_{Q}\text{-update}:&~~\lambda_{Q}^{k+1}-\lambda_{Q}^k-\rho(Ax_{[Q]}^{k+1}+Bz_{Q}^{k+1})=0\subeqn\label{eqn:af}.
\end{align}

With this formulation, we can use similar approaches in \cite{Wei:ADMMinDCO, Mateos:dissparselr} to show that (\ref{eqn:admmupdatesQ}) leads to  (\ref{eqn:distributedversionQ}) if $\lambda_Q^0$ and $z_Q^0$ are properly initialized. First multiplying the two sides of (\ref{eqn:af}) by $B^T$ and adding it to (\ref{eqn:ae}), we have
%Now we also initialize $\beta^0$ to lie in the column space of $M_-^T$ (e.g., $\beta^0=0$) such that $\beta^{k+1}$ also lies in the column space of $M_-^T$. Substitute $z^k=\frac{M_+^Tx^k}$ into (\ref{eqn:ca}) and (\ref{eqn:cb}) and we have 
\begin{align}
B^T\lambda^{k+1}_{Q} &= 0,\refstepcounter{equation}\subeqn\label{eqn:da}\\
B^T\lambda_{Q}^k+\rho B^T(Ax_{[Q]}^{k+1}+Bz_{Q}^{k+1}) &= 0.\subeqn\label{eqn:db}
\end{align}
Recalling $\lambda_{Q} =[\beta_{Q};\gamma_{Q}]$ and $B=[-I_{2EM};-I_{2EM}]$, if we initialize $\beta_Q^0=-\gamma_Q^0$, then $\beta_Q^k=-\gamma_Q^k$ for $k = 0,1,\cdots$, and (\ref{eqn:db}) implies $z_{Q}^{k+1}=\frac{1}{2}M_+^Tx_{[Q]}^{k+1}$. By initializing $z^0_{Q}=\frac{1}{2}M_+^Tx_{[Q]}^{0}$, we have $z^{k}_{Q}=\frac{1}{2}M_+^Tx_{[Q]}^{k}$ for $k = 0, 1,\cdots$. To summarize, with the initialization $\beta^0_Q=-\gamma^0_Q$ and $z^0_{Q}=\frac{1}{2}M_+^Tx_{[Q]}^{0}$, (\ref{eqn:da}) and (\ref{eqn:db}) are equivalent to 
\begin{align}
B^T\lambda^{k}_{Q} &= 0,\refstepcounter{equation}\subeqn\label{eqn:ea}\\
z^k_{Q}-\frac{1}{2}M_+^Tx^k_{[Q]} &= 0.\subeqn\label{eqn:eb}
\end{align}

%The purpose here is to study the $x$-update through the ideal ADMM which makes the problem easier.
Next we consider updates (\ref{eqn:aa})-(\ref{eqn:ac}). Multiplying the two sides of the $\lambda$-update by $A^T$ and $B^T$ and adding them to the $x$-update and $z$-update respectively, we get $\partial f(x^{k+1})+A^T\lambda^{k+1}+\rho A^TB(z^k_{Q}-z^{k+1})=0$ and $B^T\lambda^{k+1}=0$. Therefore, (\ref{eqn:aa})-(\ref{eqn:ac}) can be equivalently expressed as
\begin{align}
\partial f(x^{k+1})+A^T\lambda^{k+1}+\rho A^TB(z_{Q}^k-z^{k+1})&=0\refstepcounter{equation}\subeqn\label{eqn:ba},\\
B^T\lambda^{k+1}&=0\subeqn\label{eqn:bb},\\
\lambda^{k+1}-\lambda_{Q}^k-\rho(Ax^{k+1}+Bz^{k+1})&=0\subeqn\label{eqn:bc}.
\end{align}
Also by the definitions of $\lambda$ and $B$, we know that $\beta^{k+1}=\gamma^{k+1}$ from (\ref{eqn:bb}), and that (\ref{eqn:ba}) splits into two equations $\beta^{k+1}-\beta^k_{Q}-\rho A_1x^{k+1}+\rho z^{k+1}=0$ and $\gamma^{k+1}-\gamma_{Q}^k-\rho A_2x^{k+1}+\rho z^{k+1}=0$. Summing and subtracting these two equations lead to $z^{k+1}-\frac{1}{2}M_+^Tx^{k+1}=0$ and $\beta^{k+1}-\beta^k_{Q}-\frac{\rho}{2}M_-^Tx^{k+1}=0$, respectively. Thus, (\ref{eqn:ba})-(\ref{eqn:bc}) reduce to
\begin{align}
\partial f(x^{k+1})+A^T\lambda^{k+1}+\rho A^TB(z_{Q}^k-z^{k+1})&=0\refstepcounter{equation}\subeqn\label{eqn:ca},\\
z^{k+1}-\frac{1}{2}M_+^Tx^{k+1}&=0,\subeqn\label{eqn:cb}\\
\beta^{k+1}-\beta^k_{Q}-\frac{\rho}{2}M_-^Tx^{k+1}&=0\subeqn\label{eqn:cc}.
\end{align}
By substituting (\ref{eqn:eb}) and (\ref{eqn:cb}) into (\ref{eqn:ca}), we have (\ref{eqn:aa})-(\ref{eqn:af}) finally equivalent to
\begin{align}
\partial f(x^{k+1})+M_-\beta^{k+1}-\frac{\rho}{2}M_+M_+^Tx_{[Q]}^k+\frac{\rho}{2}M_+M_+^Tx^{k+1}=0&,\refstepcounter{equation}\subeqn\label{eqn:fa}\\
\beta^{k+1}-\beta^k_{Q}-\frac{\rho}{2}M_-^Tx^{k+1}=0&,\subeqn\label{eqn:fb}\\
%x^{k+1}_{Q}=Q(x^{k+1})&,\subeqn\label{eqn:fc}\\
\beta^{k+1}_{Q}-\beta_{Q}^k-\frac{\rho}{2}M_-^Tx_{[Q]}^{k+1}=0&.\subeqn\label{eqn:fe}
\end{align}
If we further multiply the two sides of (\ref{eqn:fb}) by $-M_-$ and add it to (\ref{eqn:fa}), we obtain 
\begin{align}
\label{eqn:xupdatefinal}
\partial f(x^{k+1})+M_-\beta^k_{Q}+\frac{\rho}{2}(M_+M_+^T+M_-M_-^T)x^{k+1}-\frac{\rho}{2}M_+M_+^Tx_{[Q]}^{k}=0.
\end{align}
From (\ref{eqn:xupdatefinal}), we see that the update of $x^{k+1}$ relies on $M_-\beta^k_Q$ instead of $\beta^k_Q$. Hence, multiplying both sides of (\ref{eqn:fe}) by $M_-$ yields\begin{align}
\label{eqn:betaupdatefinal}
M_-\beta^{k+1}_Q-M_-\beta^k_Q-\frac{\rho}{2}M_-M_-^Tx^{k+1}_{[Q]}=0.
\end{align}
 Letting $\alpha_Q^k=M_-\beta_{Q}^k$ and recalling the definitions of $L_-$, $L_+$ and $W$, we have (\ref{eqn:xupdatefinal}) and (\ref{eqn:betaupdatefinal}) equivalent to
%\begin{equation}
%\begin{aligned}
%\nabla f(x^{k+1})+\alpha^k+2\rho W x^{k+1}-\rho L_+x^k_{[Q]}=0,&\nn\\
%x^{k+1}_{[Q]}=Q(x^{k+1}),&\nn\\
%\alpha^{k+1}-\alpha^k-\rho L_-x_{[Q]}^{k+1}=0.&\nn
%\end{aligned}
%\end{equation}
%or equivalently, 
\begin{equation}
\begin{aligned}
\label{eqn:compform}
x{\text{-update}}:~&\partial f(x^{k+1})+2\rho Wx^{k+1}+\alpha_Q^k-\rho L_+x^k_{[Q]}=0,\\
%\label{eqn:compform2}
\alpha_Q{\text{-update}}:~&\alpha_Q^{k+1}-\alpha^k_Q-\rho L_-x_{[Q]}^{k+1}=0,
\end{aligned}
\end{equation}
which is exactly the matrix form of (\ref{eqn:distributedversionQ}).

\subsection{Convergence Results: Smooth Objective Functions}
\label{sec:convergence}
For ease of presentation, we first study a simple case where objective functions only contain the smooth components, i.e., $f_i=g_i$. Then $\partial f_i = \nabla g_i$ under Assumptions \ref{asm:convex} and \ref{asm:StrLipconvex}.

We consider the effect of the rounding quantization by writing $$x^{k+1}_{[Q]}=x^{k+1}+e^{k+1}.$$ Then the $\alpha_Q$-update in the QC-ADMM iteration (\ref{eqn:compform}) is equivalent to $$\alpha_Q^{k+1}=\alpha_Q^k+\rho L_-x^{k+1}+\rho L_-e^{k+1}.$$ Compared with the C-ADMM, we see that $[x^{k+1}_{[Q]}; \alpha^{k+1}_Q]$ is obtained by performing the C-ADMM update on $[x^k_{[Q]};\alpha^k_Q]$ followed by adding an error term $[e^{k+1};\rho L_-e^{k+1}]$ which is caused by the quantization operation. Due to the possible nonlinearity of $\nabla g_i$, we cannot easily write $[x^{k+1}_{[Q]};\alpha_Q^{k+1}]$ as the sum of the $(k+1)$th C-ADMM update plus an accumulated error term. Instead, we utilize the vector $u^k_Q=[z^k_Q;\beta^k_Q]$ to study the convergence of the QC-ADMM.

From Section \ref{sec:connection}, we know that the QC-ADMM can be obtained from (\ref{eqn:aa})-(\ref{eqn:af}) by initializing $\beta_Q^0=-\gamma_Q^0$ and $z_Q^0=\frac{1}{2}M_+^Tx_{[Q]}^0$. This initialization then leads to $\beta_Q^k=-\gamma_Q^k$, $z_Q^k=\frac{1}{2}M_+^Tx_{[Q]}^k$, $\beta^{k+1}=-\gamma^{k+1}$, and $z^{k+1}=\frac{1}{2}M_+^Tx^k$ for $k=0,1,\cdots$. Therefore, the $(k+1)$th update using (\ref{eqn:aa})-(\ref{eqn:af}) is only based on $[z_Q^k;\beta_Q^k]$. Similar to studying the C-ADMM through the correspondingly initialized ADMM, we investigate the updats of $[z_Q^{k+1};\beta_Q^{k+1}]$ through (\ref{eqn:aa})-(\ref{eqn:af}) that are also correspondingly initialized. To this end, we notice the following relation between $z_Q^{k+1}$ and $z^{k+1}$: $$z_Q^{k+1}=\frac{1}{2}M_+^Tx_{[Q]}^{k+1}=\frac{1}{2}M_+^Tx^{k+1}+\frac{1}{2}M_+^Te^{k+1}=z^{k+1}+\frac{1}{2}M_+^Te^{k+1}.$$ Combining (\ref{eqn:fb}) and (\ref{eqn:fe}), we also obtain $$\beta^{k+1}_Q=\beta^{k+1}+\frac{\rho}{2}M_-^Te^{k+1}.$$ Since $[z^{k+1};\beta^{k+1}]$ is obtained by performing a standard ADMM update on $[z^k_Q;\beta^k_Q]$ as seen from (\ref{eqn:aa})-(\ref{eqn:ac}), we can represent the update of $u_Q^{k+1}$ from $u_Q^k$ as
%\begin{align}
%\label{eqn:compactformQCADMM}
 %\left[\begin{array}{c} x^{k+1} \\\alpha_Q^{k+1}\end{array} \right] =\left[\begin{matrix} \nabla f + 2\rho W & 0\\ -\rho L_- & I \end{matrix}\right ] ^{-1} \left(\left[\begin{matrix} \rho L_+ & -I \\ 0 & I\end{matrix}\right]\left(\left[\begin{array}{c} x^{k} \\\alpha^{k}\end{array} \right]+\left[\begin{array}{c} e^{k} \\ 0 \end{array} \right]\right)\right),
%\end{align} 
%Recalling the consensus ADMM update, we know $[x^{k+1};\alpha^{k+1}]$ is the result of the ideal update from $[x^k_{[Q]};\alpha^k]$ plus an error term $[0;\rho L_-e^{k+1}]$. However, different from our previous work, we cannot explicitly write due the nonlinearity of $\nabla f$. We, on the other hand, handle the vector $u^k_Q=[z^k_Q;\beta^k_Q]$.
%Since $z_Q=\frac{1}{2}M_+^Tx_Q$ and Theorem \ref{thm:aberration}, we have  we have  
%\begin{align}
%\label{eqn:QCADMMzb}
 %\left[\begin{array}{c} z_Q^{k+1} \\\beta_Q^{k+1}\end{array} \right] = h\left(\left[\begin{array}{c} z_Q^{k} \\\beta_Q^{k}\end{array} \right]\right)+\left[\begin{array}{c} \frac{1}{2}M_+^Te^{k+1} \\\frac{1}{2}M_-^Te^{k+1}\end{array} \right],
%\end{align} 
%or equivalently, 
\begin{align}
\label{eqn:QCADMMzbu}
u_Q^{k+1} = \psi(u^k_Q)+u_e^k,
\end{align}
where $u_e^k = [\frac{1}{2}M_+^Te^{k+1};\frac{\rho}{2}M_-^Te^{k+1}]$, and $\psi$ denotes the standard ADMM update as defined by (\ref{eqn:uh}). We will use this relation to write $u_Q^{k+1}$ as the sum of the $(k+1)$th ADMM update from $u_Q^0$ plus an accumulated error term caused by quantization. If the QC-ADMM starts with $\alpha_Q^0$ which is in the column space $L_-$, then the $\alpha_Q$-update implies that $\alpha_Q^k$ lies in the column space of $L_-$ for $k=0,1,\cdots$. Therefore, the corresponding ADMM update possesses the linear convergence rate [cf. Equation (\ref{eqn:uhlc})] as discussed in Section \ref{sec:lcCADMM}. Utilizing this property we are able to establish the absolute convergence and hence the convergence of the accumulated error term. We first state the boundedness of $u_e^k$ and $u_Q^k$ in the following lemma.
\begin{lemma}
\label{lem:bd}
Let $f_i=g_i$ and suppose that Assumptions \ref{asm:convex} and \ref{asm:StrLipconvex} hold. Consider the QC-ADMM algorithm. Let $\beta_Q^0$ be the unique vector in the column space of $M_-$ such that $\alpha_Q^0=M_-^T\beta_Q^0$. Let also $z_Q^0=\frac{1}{2}M_+^Tx_{[Q]}^0$ and $\gamma_Q^0=-\beta_Q^0$ in (\ref{eqn:aa})-(\ref{eqn:af}).
Then for $k=0,1,\cdots$,
\begin{align}
\label{eqn:ueukbd}
\|u_e^k\|_G\leq \tau_0~\text{and}~\|u_Q^k-u^*\|_G\leq\|u_Q^0-u^*\|_G+\left(1+\frac{1}{\eta}\right)\tau_0,
\end{align}
where $u^*=[{\bf1}_{2E}\tilde{x}^*;\beta^*]$ with ${\bf1}_{2E}=1_{2E}\otimes I_M$ and $\beta^*$ being the unique vector in the column space of $M_-$ such that $M_-^T\beta^*=\alpha^*=\nabla g({\bf1}_N\tilde{x}^*)$, $\tau_0=\frac{1}{4}\Delta\sqrt{\rho M\left(\sigma_{\max}^2(M_+)+\sigma_{\max}^2(M_-)\right)}$, and $\eta$ is defined in Lemma \ref{lem:linearconvergence}.
\end{lemma}
\begin{IEEEproof}
The boundedness of $u_e^k$ follows directly from the boundedness of $e^{k+1}$, i.e., 
\begin{align}
\label{eqn:uebddd}
\|u_e^k\|_G^2 &= \rho\left\|\frac{1}{2}M_+^Te^{k+1}\right\|_2^2+\frac{1}{\rho}\left\|\frac{1}{2}\rho M_-^Te^{k+1}\right\|_2^2\nn\\
&\leq\frac{1}{4}\rho\left(\sigma_{\max}^2(M_+)+\sigma_{\max}^2(M_-)\right)\|e^{k+1}\|_2^2\nn\\
&\leq\frac{1}{16}\rho M\Delta^2\left(\sigma_{\max}^2(M_+)+\sigma_{\max}^2(M_-)\right).\nn\\
&=\tau_0^2. 
\end{align}
Since $\alpha_Q^0$ is initialized in the column space of $L_-$, we see that $\alpha_Q^k$ also lies in the column space of $L_-$ from the $\alpha_Q$-update of the QC-ADMM. Then using (\ref{eqn:Ulinearconvergence}) and (\ref{eqn:uhlc}), we have 
\begin{align}
\|u_Q^{k+1}-u^*\|_G&=\|\psi(u_Q^k)+u_e^k-u^*\|_G\nn\\
&\leq \|\psi(u^k_Q)-u^*\|_G+\|u_e^k\|_G\nn\\
&\leq \frac{1}{1+\eta}\|u^k_Q-u^*\|_G+\|u_e^k\|_G\nn\\
&~~~~\cdots\nn\\
&\leq\frac{1}{(1+\eta)^{k+1}}\|u_Q^0-u^*\|_G+\sum_{j=0}^{k}\frac{1}{(1+\eta)^j}\|u_e^j\|_G.\nn
\end{align}
Therefore, for $k = 0,1,\cdots$, we have 
\begin{align}
\label{eqn:uQbd}
\|u_Q^k-u^*\|_G \leq \|u_Q^0-u^*\|_G+\sum_{i=0}^{\infty}\frac{1}{(1+\eta)^i}\tau_0=\|u_Q^0-u^*\|_G+\left(1+\frac{1}{\eta}\right)\tau_0.
\end{align}
\end{IEEEproof}

With this lemma, we are ready to establish our main theorem as follows.
\begin{theorem}
\label{thm:convergenceQCADMM}
Let $f_i=g_i$. Consider the QC-ADMM iteration (\ref{eqn:distributedversionQ}), and suppose that Assumptions \ref{asm:convex} and \ref{asm:StrLipconvex} hold. Initializing $\alpha_Q^0$ in the column space of $L_-$, we have
\begin{enumerate}
\item{\it Convergence:} the sequence $(x_{[Q]}^k, \alpha_{Q}^k)$ generated by (\ref{eqn:distributedversionQ}) converges to a finite value $({\bf 1}_N\tilde{x}_{Q}^*,\alpha_Q^*)$ as $k\to\infty$, where $\tilde{x}^*_{Q}$ is some vector in $\Lambda^M$ and not necessarily equal to $\tilde{x}^*_{[Q]}$. % (Note that the consensus is also reached.)
\item{\it Consensus error:} an upper bound for the consensus error is given by
		$$\|\tilde{x}^*_{Q}-\tilde{x}^*\|_2\leq \left(\frac{1}{2}+\rho\frac{2E}{\sum_{i=1}^N m_{g_i}}\right)\sqrt{M}\Delta.$$
\item{\it Number of iterations:} $(x_{[Q]}^k, \alpha_{Q}^k)$ converges within $\lceil\log_{1+\eta}{\Omega}\rceil$ iterations, where $$\Omega = \max\left\{\frac{3\sqrt{\rho}\sigma_{\max}(M_-)(1+\eta)^2\left(\|u_Q^0-u^*\|_G+\tau_0\right)}{\eta\Delta},\frac{3(1+\eta)\left(\|u_Q^0-u^*\|_G+\tau_0\right)}{\sqrt{2\rho E}\eta\Delta}\right\},$$ and $\lceil y \rceil, y\in\mathbb{R},$ means the smallest integer that is greater than or equal to $y$.,
\end{enumerate}
\end{theorem}
\begin{IEEEproof} We prove the three claims one by one.
	
{\it Convergence:~}From (\ref{eqn:aa})-(\ref{eqn:af}), we know that $x^{k+1}$ is updated only based on $u_Q^k=[z_Q^k;\beta^k_Q]$ if the updates are initialized with $\beta^0_Q=-\gamma^0_Q$ and $z^0_Q=\frac{1}{2}M_-^Tx^0_{[Q]}$. As long as $u_Q^k$ converges, we must have $x_{[Q]}^{k+1}$ and $\alpha_Q^{k+1}$ converging. Given that $\alpha_Q^k$ converges, $L_-x_Q^{k}$ must converge to $0$ due to the $\alpha_Q$-update, and hence $x_{[Q]}^k$ converges to a consensus by Lemma \ref{lem:LminusP}. Therefore, to prove the convergence of $(x_{[Q]}^k,\alpha^k_Q)$, it is enough to show the convergence of $u_Q^k$.

Following (\ref{eqn:QCADMMzbu}), we have 
\begin{align}
\label{eqn:uQe}
u_Q^{k+1}&=\psi(u_Q^k) + u_e^k\nn\\
&=\psi\left (\psi(u_Q^{k-1})+u_e^{k-1}\right)+u_{e'}^k \nn\\
&=\psi^2(u_Q^{k-1}) + u_{e'}^{k-1} + u_{e'}^k\nn\\
%&=g^2\left(g(\bm u^{k-2}) +\bm u_e^{k-2}\right) + \bm u_{e'}^{k-1} + \bm u_{e'}^k\nn\\
&=\cdots\nn\\
&=\psi^{k+1}(u_Q^0) + \sum_{i=0}^k u_{e'}^{k-i},
\end{align}
where $\psi^i(\cdot), i=0,1,\cdots,$ denotes the $i$-th standard ADMM update on its argument and $u_{e'}^{k-i}=\psi^i\left(\psi(u^{k-i})+u_e^{k-i}\right)-\psi^{i+1}(u^{k-i})$. We only need to prove the convergence of the accumulated error term $\sum_{i=0}^k u_{e'}^{k-i}$ as the first term is the $(k+1)$th standard ADMM update which converges to $u^*$ as $k\to\infty$. It then suffices to show the boundedness of $\lim_{k\to\infty}\sum_{i=0}^k \|u_{e'}^{k-i}\|_G$ due to the comparison theorem and the fact that absolute convergence implies convergence \cite{Rudin:analysis}. We first obtain an upper bound on $\|u_{e'}^{k-i}\|_G$:
\begin{align}
\label{eqn:uePbd}
\|u_{e'}^{k-i}\|_{G} &= \left\|\psi^i\left(\psi(u_Q^{k-i})+u_e^{k-i}\right)-\psi^{i+1}(u_Q^{k-i})\right\|_{G}\nn\\
&= \left\|\psi^i\left(\psi(u_Q^{k-i})+ u_e^{k-i}\right)-u^*-\left(\psi^{i+1}(u_Q^{k-i})-u^*\right)\right\|_{G}\nn\\
&\leq \left\|\psi^i\left(\psi(u_Q^{k-i})+u_e^{k-i}\right)-u^*\right\|_{G}+\left\|\psi^{i+1}(u_Q^{k-i})-u^*\right\|_{G}\nn\\
&\stackrel{(a)}{\leq}\frac{1}{(1+\eta)^i}\|\psi(u_Q^{k-i})+u_e^{k-i}-u^*\|_{G}+\frac{1}{(1+\eta)^{i+1}}\|u_Q^{k-i}-u^*\|_{G}\nn\\
&\leq\frac{1}{(1+\eta)^i}\|\psi(u_Q^{k-i})-u^*\|_{G}+\frac{1}{(1+\eta)^i}\| u_e^{k-i}\|_{G}+\frac{1}{(1+\eta)^{i+1}}\|u_Q^{k-i}-u^*\|_{G}\nn\\
%&\leq\frac{2}{(1+\delta)^{i+1}}\|\bm u^{k-i}-\bm u^*\|_{\bm G}+\frac{1}{(1+\delta)^i}\|\bm u_e^{k-i}\|_{\bm G}
&\stackrel{(b)}{\leq}\frac{2}{(1+\eta)^{i+1}}\|u^0_Q-u^*\|_G + \frac{3}{(1+\eta)^i}\tau_0,
\end{align}
where ($a$) is from Lemma \ref{lem:linearconvergence} and ($b$) is due to Lemma \ref{lem:bd}. Therefore,  
\begin{align}
\label{eqn:bdaeerorterm}
 \lim_{k\to\infty}\sum_{i=0}^k \|u_{e'}^{k-i}\|_G\leq \sum_{i=0}^\infty \frac{2}{(1+\eta)^{i+1}}\|u_Q^0-u^*\|_G+\sum_{i=0}^\infty \frac{3}{(1+\eta)^i}\tau_0.
\end{align}
The convergence proof is complete by noting that $\eta >0$.

{\it Consensus error:~}The consensus error may be studied directly by calculating the accumulated error term in (\ref{eqn:uQe}). However, the bound in (\ref{eqn:bdaeerorterm}) is quite loose in general as the bounds in Lemmas \ref{lem:linearconvergence} and \ref{lem:bd} are themselves loose for the respective quantities. We alternatively study the QC-ADMM iteration (\ref{eqn:distributedversionQ}) using the fact that $x^k_{[Q]}$ converges to a consensus as $k\to\infty$.
	
Let $\tilde{x}_{Q}^*\in\Lambda^M$ be the convergent quantized value at each agent. Then $x_{i[Q]}^\infty=\tilde{x}_{Q}^*$ and $x_i^\infty=\tilde{x}_{Q}^*-e_{i}^*$ for $i = 1,2,\cdots,N$, where $e_i^*$ is the quantization error at agent $i$. It is important to note that $\tilde{x}_Q^*$ does not represent the quantized value of the global optimum $\tilde{x}^*$, i.e., $\tilde{x}_Q^*$ is not necessarily equal to $\tilde{x}^*_{[Q]}$. Summing up both sides of (\ref{eqn:distributedversionQ}) from $i = 1$ to $N$, we have 
	\begin{align}
	\label{eqn:sumbs}
	\sum_{i=1}^N \left(\nabla g_i(\tilde{x}_{Q}^*-e^*_{i})+2\rho|\mathcal{N}_i|(\tilde{x}_{Q}^*-e_{i}^*)\right) = \sum_{i=1}^N\left( |\mathcal{N}_i| \tilde{x}_{Q}^*+\rho\sum_{j\in\mathcal{N}_i} \tilde{x}_{Q}^*\right).
	\end{align}
Here we use the fact that $\alpha_Q^k$ lies in the column space of $L_-$, i.e.,  $\alpha_Q^k=L_-b^k$ for some $b^k\in\mathbb{R}^{NM}$. By Lemma \ref{lem:LminusP}, we have $$\sum_{i=1}^N \alpha^k_{iQ} = (L_-b^k)^T{\bf 1}_N=(b^k)^T(L_-^T{\bf 1}_N)=0.$$
%Denote by $x_I^*$ the solution to. Under Assumption, we have $$\sum_{i=1}^N \nabla f_i(x^*) = 0.$$ 

Since $f_i=g_i$ which is differentiable and strongly convex in $\mathbb{R}^M$, $\tilde{x}^*$ is the solution to problem (\ref{eqn:DistrOpt}) if and only if
$$\sum_{i=1}^N \nabla g_i(\tilde{x}^*) = 0.$$
Thus, (\ref{eqn:sumbs}) leads to
\begin{align}
\label{eqn:minusforlc}
\sum_{i=1}^N \left (\nabla g_i(\tilde{x}^*_{Q}-e_{i}^*) - \nabla g_i(\tilde{x}^*)\right ) = \sum_{i=1}^N 2\rho|\mathcal{N}_i|e^*_{i}.
\end{align}
Recalling the strong convexity assumption, we have
$$\|\nabla g_i(\tilde{x}^*_{Q}-e^*_{i}) - \nabla g_i(\tilde{x}^*)\|_2 \geq m_{g_i} \|\tilde{x}^*_{Q}-e_i^*-\tilde{x}^*\|_2.$$
Together with (\ref{eqn:qerr}), we obtain the following upper bound 
\begin{align}
\label{eqn:upperbd}
\|\tilde{x}^*_{Q}-\tilde{x}^*\|_2 \leq \left(\frac{1}{2}+\rho\frac{2E}{\sum_{i=1}^N m_{g_i}}\right) \sqrt{M}\Delta.
\end{align}

We next use an example to show that this bound is indeed tight. Consider a simple two-node network where $f_1(\tilde{x})=\left(\tilde{x}+\frac{3}{2}\right)^2$ and $f_2(\tilde{x})=\left(\tilde{x}+\frac{7}{2}\right)^2,\tilde{x}\in\mathbb{R}$. Then $m_{f_1}=m_{f_2}=1$ and $\tilde{x}^*=-\frac{5}{2}$. Set both $\Delta$ and $\rho$ to be $1$. In this case, we have $M=1$, $N=2$, $E=1$ and $$L_-=\begin{bmatrix} 1& -1\\ -1& 1\end{bmatrix}.$$ We start with $x^0_{1[Q]}=x^0_{2[Q]}=-1$ and $\alpha_{1Q}^0= -\alpha_{2Q}^0=1$. One can easily check that $\alpha_Q^0=[\alpha_{1Q}^0;\alpha_{2Q}^0]$ lies in the column space of $L_-$, and that $x^k_{1[Q]}= x^k_{2[Q]}=-1$ and $\alpha_{1Q}^k=-\alpha_{2Q}^k=1$ for $k =0,1,\cdots$, in the updates of (\ref{eqn:distributedversionQ}). Hence the consensus error is 
	\begin{align}
	\left \|\tilde{x}_{Q}^*-\tilde{x}^*\right\|_2=\frac{3}{2}=\left(\frac{1}{2}+\rho\frac{2E}{\sum_{i=1}^N m_{f_i}}\right)\sqrt{M}\Delta,\nn
	\end{align}
	
{\it Number of iterations:}
	The convergence of $(x_{[Q]}^k,\alpha_Q^k)$ implies that there exists a finite $k_0$ such that $\|\alpha_Q^{k}-\alpha_Q^*\|_2< \Delta$ and $\|x^{k}_{[Q]}-{\bf 1}_N\tilde{x}^*_{Q}\|_2< \Delta$ for $k\geq k_0$. Then $\alpha_Q^k=\alpha_Q^*$ and $x^k_{[Q]}={\bf 1}_N\tilde{x}^*_{Q}$ due to the rounding quantization scheme. Therefore, $(x_{[Q]}^k,\alpha_Q^k)$ converges in finite iterations.

	%h^{k+1}(u_Q^0) + \sum_{i=0}^k u_{e'}^{k-i},
	
	One may have noticed that the two terms in (\ref{eqn:uQe}) converge relatively fast: $\psi^{k+1}(u_Q^0)$ is the $(k+1)$th update of the ADMM and thus linearly convergent; $\sum_{i=0}^k u_{e'}^{k-i}$ is absolutely bounded by the sum of two geometric series  whose common ratios are both positive and less than $1$. As such, we expect an upper bound for the number of iterations that guarantees the convergence of $(x_{[Q]}^k,\alpha_Q^k)$ to $({\bf 1}_N\tilde{x}_{Q}^*,\alpha_Q^*)$. %Our purpose here is to give a rough idea of how fast the convergence can be; we will not make specific efforts to tighten the bound.  
	
	We first consider the number of iterations, denoted by $k_1$, that guarantees the convergence of $\alpha_Q^k$. Write $\beta_Q^k=\beta^k_{IQ}+\beta^k_{EQ}$ where $\beta^k_{IQ}$ and $\beta^k_{EQ}$ are the corresponding vectors in the standard ADMM update $\psi^{k}(u_Q^0)$ and the accumulated error term $\sum_{i=0}^k u_{e'}^{k-i}$ for $k\geq 1$, respectively. Define $u_Q^k=u_{IQ}^k+u_{EQ}^k$ analogously. Then we have
	\begin{align}
	\|\alpha_Q^k-\alpha_Q^*\|_2 &= \|M_-(\beta_{Q}^k-\beta_{Q}^*)\|_2\nn\\
	&=\|M_-(\beta^k_{IQ}-\beta_{IQ}^*+\beta^k_{EQ}-\beta_{EQ}^*)\|_2\nn\\
	&\leq \sigma_{\max}(M_-)\left(\|\beta^k_{IQ}-\beta_{IQ}^*\|_2+\|\beta^k_{EQ}-\beta_{EQ}^*\|_2\right)\nn\\
%	&\leq  \|\alpha_{I}^k-\alpha_I^*\|_2+\|\alpha_{E}^k-\alpha_E^*\|_2\nn\\
	\label{eqn:upperbound1}
	&\leq \sqrt{\rho}{\sigma_{\max}}(M_-)\left(\|u^k_{IQ}-u^*\|_G+\lim_{\kappa\to\infty}\|\sum_{i=k+1}^{\kappa} u_{e'}^{\kappa-i}\|_G\right),
	\end{align}
where the last inequality is from the definition of $G$. Since $\alpha_Q^k$ lies in the column space of $L_-$, we have from Lemma \ref{lem:linearconvergence} that
\begin{align}
\label{eqn:bdd1}
\|u_{IQ}^k-u^*\|_G\leq{\frac{1}{(1+\eta)^k}}\|u_Q^0-u^*\|_G.
\end{align}
Using the upper bound of (\ref{eqn:uePbd}), we get 
\begin{align}
\lim_{\kappa\to\infty}\|\sum_{i=k+1}^{\kappa} u_{e'}^{\kappa-i}\|_2&\leq \lim_{\kappa\to\infty}\sum_{i=k+1}^{\kappa} \|u_{e'}^{\kappa-i}\|_2\nn\\
&\leq \sum_{i=k+1}^\infty \left(\frac{2}{(1+\eta)^{i+1}}\|u_Q^0-u^*\|_G+\frac{3}{(1+\eta)^{i}}\tau_0\right)\nn\\
\label{eqn:bdd2}
&\leq \sum_{i=k}^\infty \left(2\|u_Q^0-u^*\|_G+3\tau_0\right)\frac{1}{(1+\eta)^{i}}.
\end{align}
Combining (\ref{eqn:bdd1}) and (\ref{eqn:bdd2}) yields
\begin{align}
\label{eqn:bdd3}
\|u^k_{IQ}-u^*\|_G+\lim_{\kappa\to\infty}\|\sum_{i=k+1}^{\kappa} u_{e'}^{\kappa-i}\|_2\leq 3\left(\|u_Q^0-u^*\|_G+\tau_0\right) \sum_{i=k}^\infty \frac{1}{(1+\eta)^{i}}.
\end{align}
Hence, it suffices to pick $k_1$ such that 
$$ 3\left(\|u_Q^0-u^*\|_G+\tau_0\right) \sum_{i=k_1}^\infty \frac{1}{(1+\eta)^{i}} < \frac{\Delta}{\sqrt{\rho}\sigma_{\max}(M_-)},$$
or, 
$$k_1= \left\lceil\log_{1+\eta}\frac{3\sqrt{\rho}\sigma_{\max}(M_-)(1+\eta)\left(\|u_Q^0-u^*\|_G+\tau_0\right)}{\eta\Delta}\right\rceil.$$

Though we have obtained $\alpha_Q^k=\alpha_Q^*$ for $k\geq k_1$, it is not enough to conclude $x_{[Q]}^{k_1}={\bf 1}_N\tilde{x}^*_{Q}$. We next use the fact that $x^{k_1+1}_{[Q]}$ reaches a consensus to find the number of iterations that guarantees the convergence of $x_{[Q]}^k$. Since $\alpha_Q^{k_1}$ has converged, we can write $x_{[Q]}^k={\bf 1}_N\zeta^k$ with $\zeta^k\in\Lambda^M$ for $k\geq k_1+1$. Thus $z_Q^k=\frac{1}{2}M_+^Tx_{[Q]}^k={\bf 1}_{2E}\zeta^k$ also reaches a consensus. Recalling that the $(k+1)$th update of (\ref{eqn:aa})-(\ref{eqn:af}) (which is properly initialized as discussed in Section \ref{sec:connection}) is only based on $[z_Q^k;\beta_Q^k]$, we only need to find $k_2$ such that $z_Q^k$ reaches ${\bf 1}_{2E}\tilde{x}_{Q}^*$ for $k\geq k_2$. Using the definitions of $u^k_Q$ and $G$, we get
\begin{align}
\|z_Q^k-z_Q^*\|_2 &\leq\frac{1}{\sqrt{\rho}}\|u_Q^{k}-u_Q^*\|_G\nn\\
&\leq \frac{1}{\sqrt{\rho}}\left(\|u^k_{IQ}-u^*\|_G+\lim_{K_1\to\infty}\|\sum_{i=k+1}^{K_1} u_{e'}^{K_1-i}\|_G\right)\nn\\
&\stackrel{(a)}\leq\frac{3}{\sqrt{\rho}}\left(\|u_Q^0-u^*\|_G+\tau_0\right)\sum_{i=k}^\infty\frac{1}{(1+\eta)^i},\nn
\end{align}
where ($a$) is from (\ref{eqn:bdd3}). Assume that $k_2\geq k_1+1$, i.e., $z_Q^k$ reaches a consensus for $k\geq k_2$. We can pick $k_2$ such that 
\begin{align}
\label{eqn:bdd4}
\|\zeta^{k_2}-\tilde{x}_{Q}^*\|=\frac{1}{\sqrt{2E}}\|z_Q^{k_2}-z_Q^*\|_2\leq \frac{3}{\sqrt{2\rho E}}\left(\|u_Q^0-u^*\|_G+\tau_0\right)\sum_{i=k_2}^\infty\frac{1}{(1+\eta)^i}<\Delta,
\end{align}
or
$$k_2= \left\lceil\log_{1+\eta}\frac{3(1+\eta)\left(\|u_Q^0-u^*\|_G+\tau_0\right)}{\sqrt{2\rho E}\eta\Delta}\right\rceil.$$
If the above $k_2$ is less than $k_1+1$, then (\ref{eqn:bdd4}) must hold by picking $k_2=k_1+1$ since the consensus of $z_Q^k$ is reached  for $k\geq k_1+1$. 

In summary, the QC-ADMM converges within $\max\{k_1+1,k_2\}$ iterations.
	%\begin{align}
	%k = \max\left\{\log_{1+\delta}\left(\frac{\rho^2 \tilde{\sigma}_{\min}^2(M_-)\sigma_{\max}^4(M_-)\|x_0-x_I^*\|_2^2+4 {\sigma}_{\max}^2(M_-)\|\alpha^0-\alpha_I^*\|_2^2}{\Delta^2 \tilde{\sigma}_{\min}^2(M_-)}\right), 2\log_{1+\delta}\left(\frac{2(1+\sqrt{\rho}\sigma_{\max}(M_-)){T_0 }}{\left(1-\sqrt{\frac{1}{1+\delta}}\right)\Delta}\right) \right\}.
	%\end{align}
	%\begin{align}
	%k = \max\left\{\log_{1+\delta}\left(\frac{4\rho\sigma_{\max}^2(M_-)\|u^k-u^0\|_G^2}{\Delta^2}\right), 2\log_{1+\delta}\left(\frac{2(1+\sqrt{\rho}\sigma_{\max}(M_-)){T_0 }}{\left(1-\sqrt{\frac{1}{1+\delta}}\right)\Delta}\right) \right\}.
	%\end{align}
	%\label{eqn:upperbound2}
	%&\leq\sqrt{\rho}{\sigma_{\max}}(M_-)\sqrt{\left(\frac{1}{1+\delta}\right)^k}\|u^0-u^*\|_G + (1+\sqrt{\rho}\sigma_{\max}(M_-))T_0\left(\frac{1}{1+\delta}\right)^k/2
	%\end{align}
\end{IEEEproof}
%As a corollary, we use an example to show that the upper bound in (\ref{eqn:upperbd}) is indeed tight.
%\begin{corollary}
%The upper bound $\left(\frac{1}{2}+\rho\frac{2E}{\sum_{i=1}^N m_{f_i}}\right) \sqrt{M}\Delta$ is tight for the consensus error.
%\end{corollary}
%\begin{IEEEproof}
%Consider a simple two-node network where $f_1(\tilde{x})=\left(\tilde{x}+\frac{3}{2}\right)^2$ and $f_2(\tilde{x})=\left(\tilde{x}+\frac{7}{2}\right)^2,\tilde{x}\in\mathbb{R}$. Then $m_{f_1}=m_{f_2}=1$ and $\tilde{x}^*=-\frac{5}{2}$. Set both $\Delta$ and $\rho$ to be $1$. In this case, we have $M=1$, $N=2$, $E=1$ and $$L_-=\begin{bmatrix} 1& -1\\ -1& 1\end{bmatrix}.$$ We start with $x^0_{1[Q]}=x^0_{2[Q]}=-1$ and $\alpha_{1Q}^0= -\alpha_{2Q}^0=1$. One can easily check that $\alpha_Q^0=[\alpha_{1Q}^0;\alpha_{2Q}^0]$ lies in the column space of $L_-$, and that $x^k_{1[Q]}= x^k_{2[Q]}=-1$ and $\alpha_{1Q}^k=-\alpha_{2Q}^k=1$ for $k =0,1,\cdots$, in the updates of (\ref{eqn:distributedversionQ}). Hence the consensus error is 
%	\begin{align}
%	\left \|\tilde{x}_{Q}^*-\tilde{x}^*\right\|_2=\frac{3}{2}=\left(\frac{1}{2}+\rho\frac{2E}{\sum_{i=1}^N m_{f_i}}\right)\sqrt{M}\Delta,\nn
%	\end{align}
%which is the error bound in (\ref{eqn:upperbd}).
%\end{IEEEproof}
\begin{remark}
\label{rmk:localoptima}
We shall mention that the limit $({\bf 1}_N\tilde{x}^*_{Q},\alpha_Q^*)$ need not be unique. This is because, unlike the standard ADMM, $\|u_Q^k-u^*\|_G$ in the QC-ADMM need not decrease monotonically due to the quantization that occurs on $x^k$ at each update. Note also that the given example illustrating the tightness of the consensus error bound is {\emph{poorly}} initialized and we usually have smaller errors than the upper bound of (\ref{eqn:upperbd}) in practice (see simulations). 
\end{remark}
\begin{remark}
An interesting observation of the above theorem is the parameters $\mu$ given in Lemma \ref{lem:linearconvergence} and $\rho$ which is the step size of the dual variable update. %We can see that $\mu$ involves in the convergence rate $\eta$ of the C-ADMM while $\rho$ affects the convergence rate $\eta$, the consensus error and also the convergence time of the QC-ADMM. 
Without the knowledge of the network topology or the properties of local objective functions (namely, $m_g$ and $M_g$ defined in Lemma \ref{lem:linearconvergence}), one can hardly determine the optimal selection of these parameters under any criteria. In this sense, we do not regard $\mu$ and $\rho$ as a factor affecting the performance of the QC-ADMM, but may simply set, e.g., $\mu=\frac{3}{2}$ and $\rho=1$. Nevertheless, we will simulate the QC-ADMM with different $\rho$ for a distributed LASSO problem in Section \ref{sec:simulations}.
\end{remark}
%\begin{remark}
%\label{rmk:rho}
%An interesting observation of our main result is the parameter $\rho$. 
%While $\rho$ directly affects the consensus error bound as seen from (\ref{eqn:upperbd}), it is not easy to characterize how it affects the convergence time. Indeed one can simply set $\rho=1$, and we do not regard $\rho$ as a factor affecting our algorithm's performance. We refer readers to  \cite{Boyd:ADMM,Wei:ADMMinDCO,Ghadimi:optimalpselctio} for discussions on the effect of $\rho$ on the ADMM. %Since $N$ and $\sigma_{\max}(M_-)$ are greater than $1$, we have $k_2=2$ which is relatively small compared with $k_1$. We refer readers to \cite{Boyd:ADMM,Wei:ADMMinDCO,Ghadimi:optimalpselctio} for how $\rho$ affects the ADMM's performance.  
%\end{remark}
%\begin{remark}
%The number of iterations that guarantees convergence is bounded by $\log_{1+\delta}\Omega$. This bound seems more favorable in large scale networks than the convergence time characterized by a polynomial of the number of nodes, quantization resolution and agents' data in \cite{Nedic:DAAQeffects}. However, $\delta$ can be very small when the graph density is high thus resulting in slow convergence rate. See simulations in Figure~\ref{fig:ctime}
%\end{remark}
\begin{remark}
\label{eqn:PfIdeaGeneralization}
To show the convergence of the QC-ADMM, our proof utilizes the linear convergence of the ADMM update $\psi$ on $u_Q^k$ and the boundedness of the error term $u_{e}^k$. As such, the main result for rounding quantization also holds for other deterministic quantizations as long as the quantization error is bounded. The proof idea also works for proving the convergence of other distributed algorithms with inexact updates: if the exactly updated variables converge relatively fast (e.g., $O(1/k^2)$ and linear rate) and the error term on these variables at each update is deterministic and bounded, then this algorithm must converge. This idea will be used in Section \ref{sec:convergenceG} to show the convergence results of the general objective functions that may contain the non-smooth components.

%as our proof only utilized the deterministic scheme and bounded quantization error.
\end{remark}

\begin{remark}
\label{rmk:otherscheme}
As previously mentioned, $x^k_i$ in the QC-ADMM is quantized for the $(k+1)$th update at its own agent even though agents can compute and store real values with infinite precision. The reason is to guarantee that $\alpha_Q^k$ lies in the column space of $L_-$ and that the QC-ADMM possesses the linear convergence rate at the ADMM update $\psi$ on $u_Q^k$ [cf. Equation (\ref{eqn:QCADMMzbu})].
\end{remark}

\subsection{Convergence Results: General Objective Functions}
\label{sec:convergenceG}
We now investigate the general case where $f_i = g_i + h_i$. Even though (\ref{eqn:QCADMMzbu}) still represents the update of $u_Q^{k+1}$ from $u_Q^k$, the proof for smooth objective functions does not apply since the ADMM update $\psi$ no longer preserves the linear convergence property [cf. Equation (\ref{eqn:Ulinearconvergence})] due to the non-smooth component $h_i$. To proceed, we write explicitly $\partial f_i = \nabla g_i + \partial h_i$. Then the $x$-update (\ref{eqn:aa}) is equivalent to
\begin{align}
\label{eqn:xupdategeneral}
%x^{k+1}=\argmin_{x} f(x)+g(x)+\langle\lambda_Q^k,Ax+Bz^k_Q\rangle+\frac{\rho}{2}\|Ax+Bz\|_2^2
\nabla g(x^{k+1})+\partial h(x^{k+1})+ A^T\lambda_Q^k+\rho A^T(Ax^{k+1}+Bz_Q^k)=0.
\end{align} 

We next make a further assumption on the smooth component $g_i$:
\begin{assumption}
\label{asm:GUsolution}
The solution to $\min_{\tilde{x}}\sum_{i=1}^N g_i(\tilde{x})$ is attainable.
\end{assumption}

This assumption together with Assumption \ref{asm:StrLipconvex} implies that there exists a unique ${u'}^*=[{z'}^*;{\beta'}^*]$ which is the optimal value of solving $\min_{\tilde{x}}\sum_{i=1}^N g_i(\tilde{x})$ using the properly initialized ADMM. Next we consider the ADMM update on $u_Q^k=[z_Q^k;\beta_Q^k]$ where only $g_i$'s are involved, i.e., 
\begin{align}
\label{eqn:xupdategeneralsmooth}
%x^{k+1}=\argmin_{x} f(x)+g(x)+\langle\lambda_Q^k,Ax+Bz^k_Q\rangle+\frac{\rho}{2}\|Ax+Bz\|_2^2
\nabla g({x'}^{k+1})+ A^T\lambda_Q^k+\rho A^T(A{x'}^{k+1}+Bz_Q^k)=0,
\end{align} 
and the subsequent updates follow from (\ref{eqn:ab})-(\ref{eqn:af}). Here we use $x'^{k+1}$, ${u'}^{k+1}$, and ${u'}_Q^{k+1}$ to denote the corresponding updated values from (\ref{eqn:xupdategeneralsmooth}) and (\ref{eqn:ab})-(\ref{eqn:af}) based on $u_Q^k=[z_Q^k;\beta_Q^k]$. Denote by $\psi'$ the update ${u'}^{k+1}$ from $u_Q^k$, i.e., ${u'}^{k+1}=\psi'(u_Q^k)$. By Theorem \ref{lem:linearconvergence}, we have 
\begin{align}
\label{eqn:psiprime}
\|{u'}^{k+1}-{u'}^*\|_G \leq \frac{1}{1+\eta}\|u_Q^k-{u'}^*\|_G.
\end{align}

To apply the previous proof idea, we rewrite (\ref{eqn:QCADMMzbu}) as 
\begin{align}
u_Q^{k+1} = \psi'(u_Q^k) + \psi(u_Q^k) - \psi'(u_Q^k) + u_e^k.\nn
\end{align}
Treating  $\psi(u_Q^k) - \psi'(u_Q^k) + u_e^k$ as the error term, we see from Lemma \ref{lem:bd} and Theorem \ref{thm:convergenceQCADMM} that the QC-ADMM must converge when $\psi(u_Q^k) - \psi'(u_Q^k) + u_e^k$ is bounded throughout the updates. Note that $u_e^k$ is bounded as a result of the boundedness of $e^{k+1}$ [cf. Lemma \ref{lem:bd}], and that
\begin{align}
\psi(u_Q^k)-\psi'(u_Q^k) = u^{k+1}-{u'}^{k+1}=\left[\begin{array}{c}\frac{1}{2}M_ +^T\\ \frac{1}{2}M_-^T \end{array} \right]\left(x^{k+1}- {x'}^{k+1}\right).\nn
\end{align}
The QC-ADMM thus converges as long as $(x^{k+1} - {x'}^{k+1})$ is bounded, which is stated below.

\begin{theorem}
\label{thm:generalcase}
Consider the QC-ADMM algorithm. Suppose that Assumptions \ref{asm:convex}--\ref{asm:GUsolution} hold. If $\|{x'}^{k+1}-x^{k+1}\|_2\leq \Delta_x$ for some $\Delta_x>0$ throughout the iterations, we have
\begin{enumerate}
\item{\it Convergence:} the sequence $(x_{[Q]}^k, \alpha_{Q}^k)$ generated by (\ref{eqn:distributedversionQ}) converges to a finite value $({\bf 1}_N\tilde{x}_{Q}^*,\alpha_Q^*)$ as $k\to\infty$, where $\tilde{x}^*_{Q}\in\Lambda^M$. % (Note that the consensus is also reached.)
\item{\it Consensus error:} an upper bound for the consensus error is given by
		$$\|\tilde{x}^*_{Q}-\tilde{x}^*\|_2\leq \left(\frac{1}{2}+\rho\frac{2E}{\sum_{i=1}^N m_{g_i}}\right)\sqrt{M}\Delta.$$
\item{\it Number of iterations:} $(x_{[Q]}^k, \alpha_{Q}^k)$ converges within $\lceil\log_{1+\eta}{\Omega}\rceil$ iterations, where  $$\tau_1=\left(\frac{1}{2}\Delta_x+\frac{1}{4}\Delta\sqrt{M}\right)\sqrt{\rho\sigma_{\max}^2(M_+)+\frac{1}{\rho}\sigma_{\max}^2(M_-)},$$ and 
\begin{align}
\label{eqn:OmegaUB}
\Omega = \max\left\{\frac{3\sqrt{\rho}\sigma_{\max}(M_-)(1+\eta)^2\left(\|u_Q^0-{u'}^*\|_G+\tau_1\right)}{\eta\Delta},\frac{3(1+\eta)\left(\|u_Q^0-{u'}^*\|_G+\tau_1\right)}{\sqrt{2\rho E}\eta\Delta}\right\}.
\end{align}
\end{enumerate}
\end{theorem}
\begin{IEEEproof}
We only outline the proof as it is similar to that of Theorem \ref{thm:convergenceQCADMM}.

{\it Convergence:} We first have 
\begin{align}
\|\psi(u_Q^k) - \psi'(u_Q^k) + u_e^k\|_G &\leq \left\|\left[\begin{array}{c}\frac{1}{2}M_ +^T\\ \frac{1}{2}M_-^T \end{array} \right]\left(x^{k+1}- {x'}^{k+1}\right)\right\|_G+\|u_e^k\|_G\nn\\
&\stackrel{(a)}{\leq} \left(\frac{1}{2}\Delta_x+\frac{1}{4}\Delta\sqrt{M}\right)\sqrt{\rho\sigma_{\max}^2(M_+)+\frac{1}{\rho}\sigma_{\max}^2(M_-)}\nn\\
&=\tau_1,\nn
\end{align}
where $(a)$ is due to the definition of $G$-norm and (\ref{eqn:uebddd}). Then one can similarly obtain an upper bound on $u_Q^k$ which is
$$\|u_Q^k\|_G\leq\|u_Q^0-u'^*\|_G+\left(1+\frac{1}{\eta}\right)\tau_1.$$
The rest of the proof is exactly the same as Theorem \ref{thm:convergenceQCADMM}.

{\it Consensus error:} Replace $\nabla g_i$ with $\partial f_i$ in (\ref{eqn:sumbs}). Note that we again have $\sum_{i=1}^N \partial f_i(\tilde{x}^*) = 0,$ and (\ref{eqn:minusforlc}) becomes
$$\sum_{i=1}^N \left (\partial f_i(\tilde{x}^*_{Q}-e_{i}^*) - \partial f_i(\tilde{x}^*)\right ) = \sum_{i=1}^N 2\rho|\mathcal{N}_i|e^*_{i}.$$
Using the convexity and strong convexity assumptions, we know that 
$$\left(\partial f_i(\tilde{x}^*_{Q}-e_{i}^*) - \partial f_i(\tilde{x}^*)\right)^T(\tilde{x}^*_{Q}-e_{i}^*-\tilde{x}^*)\geq m_{g_i}\|\tilde{x}^*_{Q}-e_{i}^*-\tilde{x}^*\|_2^2,$$
and thus $$\|\partial f_i(\tilde{x}^*_{Q}-e_{i}^*) - \partial f_i(\tilde{x}^*)\|_2\geq m_{g_i}\|\tilde{x}^*_{Q}-e_{i}^*-\tilde{x}^*\|_2.$$
Then the upper bound (\ref{eqn:upperbd}) also holds.

{\it Number of iterations:} The bound for number of iterations can be obtained by replacing $u^*$ and $\tau_0$ with ${u'}^*$ and $\tau_1$, respectively.
\end{IEEEproof}

We then provide two often used non-smooth functions that satisfy the condition $\|x'^{k+1}-x^{k+1}\|_2\leq\Delta_x$: $\ell_1$-norm and indicator function with bounded box set.

{\it{$\ell_1$-norm:}} Let $\|w\|_1$ denote the $\ell_1$-norm of a vector $w\in\mathbb{R}^L$. Then its subgradient is given by
\begin{align}
%\label{eqn:subl1norm}
\left(\partial \|w\|_1\right)_j = \begin{cases}
1, &w_j > 0,\\
-1~\text{or}~1, &w_j = 0,\\
-1, &w_j < 0.
\end{cases}\nn
\end{align}
where $w_j$ is the $j$th entry of $w$. 
Consider the non-smooth component $h_i(\tilde{x}) = \xi_i\|\tilde{x}\|_1$ and hence $h(x) = \sum_{i=1}^N \xi_i\|x_i\|_1$, where $\xi_i>0$. Subtracting (\ref{eqn:xupdategeneralsmooth}) from (\ref{eqn:xupdategeneral}) we get
\begin{align}
\label{eqn:l1norm}
\nabla g(x^{k+1}) + \rho A^TAx^{k+1} - \nabla g({x'}^{k+1}) - \rho A^TA{x'}^{k+1} = -\partial h({x}^{k+1}).
\end{align}
Since $A^TA=2\rho W$, we have 
\begin{align}
\label{eqn:gx}
&\left(\nabla g(x^{k+1})-\nabla g({x'}^{k+1})\right)^T(x^{k+1}-{x'}^{k+1}) + \rho\left(A^TA(x^{k+1}-{x'}^{k+1})\right)^T(x^{k+1}-{x'}^{k+1}) \nn\\
\geq &\left(m_{g}+ 2\rho|\mathcal{N}|_{\min}\right)\|x^{k+1}-{x'}^{k+1}\|_2^2,
\end{align}
where the last inequality is due to the strong convexity of $g_i$ and $|\mathcal{N}|_{\min}=\min_i |\mathcal{N}_i|$.
Also, using the Cauchy-Schwarz inequality we get 
\begin{align}
\label{eqn:hx}
\left(-\partial h({x}^{k+1})\right)^T\left(x^{k+1}-{x'}^{k+1}\right) \leq \left\|\partial \|{x}^{k+1}\|_1\right\|_2 \|x^{k+1}-{x'}^{k+1}\|_2\leq\left({M\sum_{i=1}^N\xi_i^2}\right)^{1/2}\|x^{k+1}-{x'}^{k+1}\|_2.
\end{align}
Combining (\ref{eqn:gx}) and (\ref{eqn:hx}), we obtain that $$\|x^{k+1}-{x'}^{k+1}\|_2\leq \frac{\left({M\sum_{i=1}^N\xi_i^2}\right)^{1/2}}{m_g+2\rho |\mathcal{N}|_{\min}}.$$
 
{\it Indicator function with bounded box set}: An indicator function is usually used when the optimization variable is subject to a constraint set. For example, if $w\in\mathcal{X}\subset\mathbb{R}^L$ for some set $\mathcal{X}$, then this can be included in the indicator function defined as
\begin{align}
I_\mathcal{X}(w) =
\begin{cases}
0, & {\text{if}}~w \in \mathcal{X},\\
\infty, & \text{otherwise}.
\end{cases}
\end{align}

We consider $f_i(\tilde{x})=g_i(\tilde{x})+I_\mathcal{X}(\tilde{x})$ where $\mathcal{X}$ is a nonempty compact box set, i.e., $\mathcal{X}=\{\tilde{x}\in\mathbb{R}^M:a\preceq\tilde{x}\preceq b\}$ with $a, b\in\mathbb{R}^M$ and $\preceq$ representing the component-wise inequality. From (\ref{eqn:distributedversionQ}) it is clear that the $(k+1)$th $x_i$-update at agent $i$ is equivalent to 
\begin{align}
\label{eqn:xupdateG}
x_i^{k+1}&=\argmin_{\tilde{x}\in\mathcal{X}} g_i(\tilde{x})+\rho|\mathcal{N}_i|\tilde{x}^T\tilde{x}-\rho\left(|\mathcal{N}_i|x_{i[Q]}^k+\sum_{j\in\mathcal{N}_i}x_{j[Q]}^k\right)^T\tilde{x} + \left(\alpha_{iQ}^k\right)^T\tilde{x}\nn\\
&= \argmin_{\tilde{x}\in\mathcal{X}} G_i^k(\tilde{x}) + \left(\alpha_{iQ}^k\right)^T\tilde{x},
\end{align}
where we define $G^k_i(\tilde{x})=g_i(\tilde{x})+\rho|\mathcal{N}_i|\tilde{x}^T\tilde{x}-\rho\left(|\mathcal{N}_i|x_{i[Q]}^k+\sum_{j\in\mathcal{N}_i}x_{j[Q]}^k\right)^T\tilde{x}$ for ease of presentation. Define $Q_0=\sup\{\|\tilde{x}\|_2+\frac{1}{2}\Delta\sqrt{M}\mid\tilde{x}\in\mathcal{X}\}$. Then $\|\tilde{x}\|_2\leq Q_0$ and $\|\tilde{x}_{[Q]}\|_2\leq Q_0$ for any $\tilde{x}\in\mathcal{X}$. Since $g_i(\tilde{x})$ is convex and differentiable, there must exist $t_i\in\mathbb{R}$ such that $\|\nabla g_i(\tilde{x})\|_2\leq t_i$ for $\tilde{x}\in\mathcal{X}$. Noting that (\ref{eqn:xupdateG}) implies $x_i^{k+1}\in\mathcal{X}$, one can easily verify that $G_i^{k+1}(\tilde{x})$ is Lipschitz continuous with parameter $t_i+4\rho|\mathcal{N}_i|Q_0$. Denote the $l$th entry of $x_{i}^k$ and $\alpha_{iQ}^k$ as $x_{i_l}^k$ and $\alpha_{{iQ}_l}^k$, respectively. We claim the following: 
\begin{align}
\label{eqn:boundx}
x_{i_l}^{k+1} = 
\begin{cases}
a_i, &\text{if}~\alpha_{iQ_l}^k>t_i+4\rho|\mathcal{N}_i|Q_0,\\
b_i, &\text{if}~\alpha_{iQ_l}^k<-(t_i+4\rho|\mathcal{N}_i|Q_0),
\end{cases}
\end{align}
for $k=1,2,\cdots$.

Assume that $x_{i_l}^{k+1} > a_i$ when $\alpha_{iQ_l}^k>t_i+4\rho|\mathcal{N}_i|Q_0$. Let $\tilde{x}^{k+1}\in\mathbb{R}^M$ denote the vector with the $l$th entry being $a_i-x_{i_l}^{k+1}$ and the rest entries being $0$. We have $x_i^{k+1}+\tilde{x}^{k+1}\in\mathcal{X}$, and  
\begin{align}
\label{enq:notmin}
&G_i^k(x_i^{k+1}+\tilde{x}^{k+1})+\left(\alpha^k_{iQ}\right)^T(x_i^{k+1}+\tilde{x}^{k+1})-G_i^k(x_i^{k+1})-\left(\alpha^k_{iQ}\right)^Tx_i^{k+1}\nn\\
\stackrel{(a)}\leq&\left(t_i+4\rho|\mathcal{N}_i|Q_0\right)(x_{i_l}^{k+1}-a_i)-a_{iQ_l}^{k+1}(x_{i_l}^{k+1}-a_i)\nn\\
=&\left(\alpha_{iQ_l}^{k+1}-(t_i+4\rho|\mathcal{N}_i|Q_0)\right)\left(a_i-x_{i_l}^{k+1}\right),
\end{align}
where $(a)$ is because $G_i^k$ is Lipschitz continuous for $k\geq 1$. Since $\alpha_{iQ_l}^{k+1}>t_i+4\rho|\mathcal{N}_i|Q_0$ and $x_{i_l}^{k+1}>a_i$, (\ref{enq:notmin}) is negative, which contradicts the fact that $x_i^{k+1}$ minimizes $G_i^k(\tilde{x}) + \left(\alpha_{iQ}^k\right)^T\tilde{x}$ over $\mathcal{X}$. That $x_{i_l}^{k+1}=b_i$ when $\alpha_{iQ_l}^k<-(t_i+4\rho|\mathcal{N}_i|Q_0)$ can be similarly shown.

With (\ref{eqn:boundx}) we demonstrate that for $k=0,1,\cdots$,
\begin{align}
\label{eqn:bdalphaiq}
\|\alpha_{iQ}^k\|_2\leq\max\left\{\|\alpha_{iQ}^0\|_2+2\rho|\mathcal{N}_i|Q_0,(t_i+6\rho|\mathcal{N}_i|Q_0)\sqrt{M}\right\}.
\end{align}  
Let $k\geq1$. If $\alpha_{iQ_l}^k>t_i+4\rho|\mathcal{N}_i|Q_0$, then $x_{i_l}^{k+1}=a_i$ as a result of (\ref{eqn:boundx}). Since $x_{j}^{k+1}\in\mathcal{X}$ for $j\in\mathcal{N}_i$ and $x_{j_l}^{k+1}\geq a_i$, $\alpha_{iQ_l}^{k+1}=\alpha_{iQ_l}^k+\rho|\mathcal{N}_i|x_{i_l[Q]}^{k+1}-\rho\sum_{j\in\mathcal{N}_i}x_{j_l[Q]}^{k+1}\leq\alpha_{iQ_l}^k$. Similarly, $\alpha_{iQ_l}^{k+1}\geq\alpha_{iQ_l}^k$ if $\alpha_{iQ_l}^k<-(t_i+4\rho|\mathcal{N}_i|Q_0)$.
If $|\alpha_{iQ_l}^k|\leq t_i+4\rho|\mathcal{N}_i|Q_0$, we have that $|\alpha_{iQ_l}^{k+1}|\leq|\alpha_{iQ_l}^k|+\rho|\mathcal{N}_i|x_{i_l[Q]}^{k+1}|+\rho\sum_{l\in\mathcal{N}_i}|x_{l_j[Q]}^{k+1}|\leq t_i+6\rho|\mathcal{N}_i|Q_0$. When $k=0$, we obtain that$\|\alpha_{iQ}^{1}\|_2=\|\alpha_{iQ}^0+\rho|\mathcal{N}_i|x_{i[Q]}^{1}-\rho\sum_{j\in\mathcal{N}_i} x_{j[Q]}^{1}\|_2\leq\|\alpha_{iQ}^0\|_2+2\rho|\mathcal{N}_i|Q_0$ as $x_i^1\in\mathcal{X}$ and $x_j^1\in\mathcal{X}$ for $j\in\mathcal{N}_i$. We finally derive (\ref{eqn:bdalphaiq}) by the definition of the Euclidean norm.

Next we use the strong convexity of $g_i$ to show that $\|x'^{k+1}-x^{k+1}\|_2$ is bounded. Noting that $x_i'^{k+1}= \argmin_{\tilde{x}} G_i^k(\tilde{x}) + \left(\alpha_{iQ}^k\right)^T\tilde{x},$ we get that $\nabla G_i^k({x_i'^{k+1}})+\alpha_{iQ}^k=0$. Thus, 
\begin{align}
\nabla G_i^k(x_i^{k+1})+\alpha_{iQ}^k &=\nabla G_i(x_i^{k+1})+\alpha_{iQ}^k -\nabla G_i^k({x_i'^{k+1}})-\alpha_{iQ}^k\nn\\
&=\nabla g_i(x_i^{k+1})+2\rho|\mathcal{N}_i|x_i^{k+1}-\nabla g_i(x_i'^{k+1})-2\rho|\mathcal{N}_i|x_i'^{k+1}.\nn
\end{align}
Since $g_i$ is strongly convex under Assumption \ref{asm:StrLipconvex}, we have 
\begin{align}
\label{eqn:ss22}
\|\nabla G_i^k(x_i^{k+1})+\alpha_{iQ}^k\|_2\geq \left(m_{g_i}+2\rho|\mathcal{N}_i|\right)\|x_i^{k+1}-x_i'^{k+1}\|_2.
\end{align}
In addition, when $k\geq 1$,
\begin{align}
\label{eqn:ss11}
\|\nabla G_i^k(x_i^{k+1})+\alpha_{iQ}^k\|_2&=\left\|\nabla g_i(x_i^{k+1}) + 2\rho|\mathcal{N}_i|x_i^{k+1}-\rho\left(|\mathcal{N}_i|x_{i[Q]}^k+\sum_{j\in\mathcal{N}_i}x_{j[Q]}^k\right)+\alpha_{iQ}^k\right\|_2\nn\\
&\le t_i+4\rho|\mathcal{N}_i|Q_0+\|\alpha_{iQ}^k\|_2.
\end{align}
Combining (\ref{eqn:bdalphaiq}), (\ref{eqn:ss22}) and (\ref{eqn:ss11}), we obtain that for $k\geq 1$, $$\|x_i^{k+1}-x_i'^{k+1}\|_2\leq\frac{\max\{\|\alpha_{iQ}^0\|_2+t_i+6\rho|\mathcal{N}_i|Q_0,(\sqrt{M}+1)t_i+(4+6\sqrt{M})\rho|\mathcal{N}_i|Q_0\}}{m_{g_i}+2\rho|\mathcal{N}_i|},$$ and hence
$$\|x^{k+1}-x'^{k+1}\|_2\leq\sum_{i=1}^N\frac{\max\{\|\alpha_{iQ}^0\|_2+t_i+6\rho|\mathcal{N}_i|Q_0,(\sqrt{M}+1)t_i+(4+6\sqrt{M})\rho|\mathcal{N}_i|Q_0\}}{m_{g_i}+2\rho|\mathcal{N}_i|}.$$

%Then $x^{k+1}$ must lie in $\mathcal{X}$ and hence be bounded. Therefore, it suffices to show that ${x'}^{k+1}$ is bounded. Since $g_i(x)$ is strongly convex under Assumption \ref{asm:StrLipconvex}, (\ref{eqn:xupdategeneralsmooth}) implies that we only need that $\lambda_Q^k$ is bounded, which is equivalent to proving that $\beta_Q^k$ is bounded as $\lambda_Q^k=[\beta_Q^k;-\beta_Q^k]$. Moreover, that $\beta^k_Q$ lies in the column space of $M_-^T$ implies $\|\alpha_Q^k\|_2=\|M_-\beta_Q^k\|_2\geq\tilde{\sigma}(M_-)\|\beta_Q^k\|_2$, and hence it is enough to show that $\alpha_{iQ}^k$ is bounded. 

%For this case, I can only show that $(x^{k+1}-{x'}^{k+1})$ is bounded when $f_i(\tilde{x})=g_i(\tilde{x})+I_\mathcal{X}(\tilde{x})$ where $\mathcal{X}$ is a nonempty compact box set, i.e., $\mathcal{X}=\{\tilde{x}\in\mathbb{R}^M:l_0\preceq\tilde{x}\preceq l_1\}$ for some $l_0,l_1\in\mathbb{R}^M$.
%Consider $h_i(\tilde{x}) = I_{\mathcal{X}}(\tilde{x})$ where $\mathcal{X}$ is a nonempty compact box set, i.e., $\mathcal{X}=\{l_0\preceq\tilde{x}\preceq l_1\}$ for $l_0,l_1\in\mathbb{R}^M$. Then it is clear ${x}^{k+1}$ is bounded and we only need to show ${x'}^{k+1}$ is bounded. Since $z_Q^{k+1} = 

%\begin{enumerate}
%\item $\ell_1$-norm
%We consider $g_i(\tilde{x})=\psi_i\|\tilde{x}\|_1$. And therefore $g(x)=We then have the subgradient given by 
%\begin{align}
%x_N\text{-update}:&~~\nabla f(x_N^{k+1}) + \nabla fA^T\lambda_{QN}^k+\rho A^T(Ax^{k+1}+Bz_{QN}^k)=0,\nn\\
%x_N\text{-update}:&~~\nabla f(x_N^{k+1}) + A^T\lambda_{QN}^k+\rho A^T(Ax^{k+1}+Bz_{QN}^k)=0,\nn
%\end{align}
%\end{enumerate}

\section{Simulations}
\label{sec:simulations}
In this section, we present some simulation results to examine previous theoretical analysis. 

\begin{figure}[htbp]
	\centering
	\includegraphics[width=0.7\textwidth]{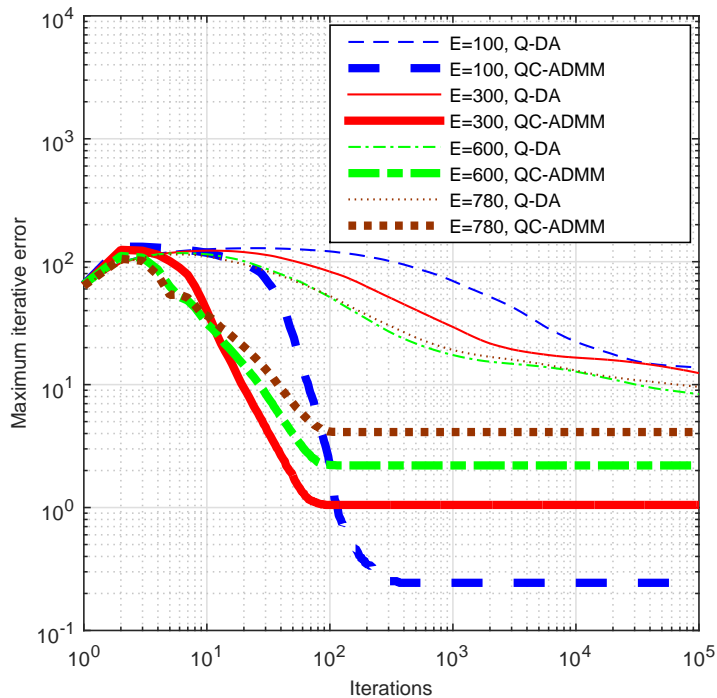}%
		\caption{Performance of the QC-ADMM and Q-DA where the plotted values are the average of $100$ runs.}
	\label{fig:DAADMM}
\end{figure}

To construct a connected network with $N$ nodes and $E$ edges, we generate a complete graph consisting of $N$ nodes, and then randomly remove $\frac{N(N-1)}{2}-E$ edges while ensuring that the network stays connected. We first consider the following distributed optimization problem:
\begin{equation}
\begin{aligned}
\label{eqn:admm}
& \underset{\tilde{x}}{\text{minimize}}
& & \sum_{i=1}^N |a_i|\|\tilde{x}\|_2^2 + b_i^T\tilde{x}\nn\\
& \text{subject to}
& & 1_M \preceq \frac{\tilde{x}}{N} \preceq 1_M,\nn\\
%& & x \in \mathbb{R}^M
\end{aligned}
\end{equation}
where $a_i\in\mathbb{R}\setminus\{0\}$ follows $\mathcal{N}(0,1)$ and $b_i\in\mathbb{R}^M$ has its entries following $\mathcal{N}(0,N^4)$. Since we do not assume $\emph{a priori}$ knowledge of the network structure, the quantized incremental algorithm does not work here. Also, the quantized subgradient method in \cite{Nedic:Dsubgradient2} does not have results for constrained problems. Hence we can only use the QC-ADMM and the quantized dual averaging (Q-DA) method. Let $\Delta=1$ and $M=3$. Set $\rho=1$ and the initial variables $x_i^0=\alpha_{iQ}^0=0$ for the QC-ADMM. The proximal function is chosen as $e(x)=\frac{1}{2}\|x\|_2^2$ for the Q-DA. Fig.~\ref{fig:DAADMM} shows the simulation results of a network with $N=40$ nodes where the maximum iterative error is defined by $\max_{i=1}^N {\|x_{i[Q]}^k-\tilde{x}^*\|_2}$.

As seen from Fig.~\ref{fig:DAADMM}, the QC-ADMM has small maximum iterative errors and converges fast. Note that the QC-ADMM converges to a consensus in finite iterations, while the Q-DA only reaches a neighborhood of $\tilde{x}^*$ and does not guarantee the convergence or a consensus (see \cite{Yuan:DdualQ}). We also check that the practical consensus error of the QC-ADMM is usually much smaller than the upper bound in Theorem \ref{thm:generalcase}. For example, when $E=300$, the average upper bound in the simulation is $18.00$ while the practical average consensus error is $1.05$. Another interesting observation is that when the graph becomes denser, i.e., $E$ becomes larger, the consensus error of the QC-ADMM tends to increase, which is in accordance with the upper bound for the consensus error.

We next consider a distributed LASSO problem: 
\begin{equation}
\begin{aligned}
\label{eqn:admm}
& \underset{\tilde{x}}{\text{minimize}}\nn
& & \sum_{i=1}^N \|A_i\tilde{x}-y_i\|_2+\lambda_i\|\tilde{x}\|_1, \nn
\end{aligned}
\end{equation}
where $A_i\in\mathbb{R}^{M\times M}$ is the linear measurement matrix of agent $i$ whose elements follow $\mathcal{N}(0,1)$, $y_i\in\mathbb{R}^{M}$ is the measurement vector of agent $i$ whose elements follow $\mathcal{N}(0,N^2)$, and $\lambda_i\in\mathbb{R^+}$ is a positive weight at agent $i$ and follows $\mathcal{N}(0,N^2)$. Let $M=20$ and $x_i^0=\alpha_{iQ}^0=0$. Define the iterative error as
$\|x_{[Q]}^k-{\bf 1}_N\tilde{x}^*\|_2/\sqrt{N}$ which is equal to the consensus error when a consensus is reached. In the following we study the effects of the quantization resolution, the algorithm parameter, and the graph density on the QC-ADMM via this LASSO problem. 

\begin{figure}[htbp]
	\centering
	\includegraphics[width=0.7\textwidth]{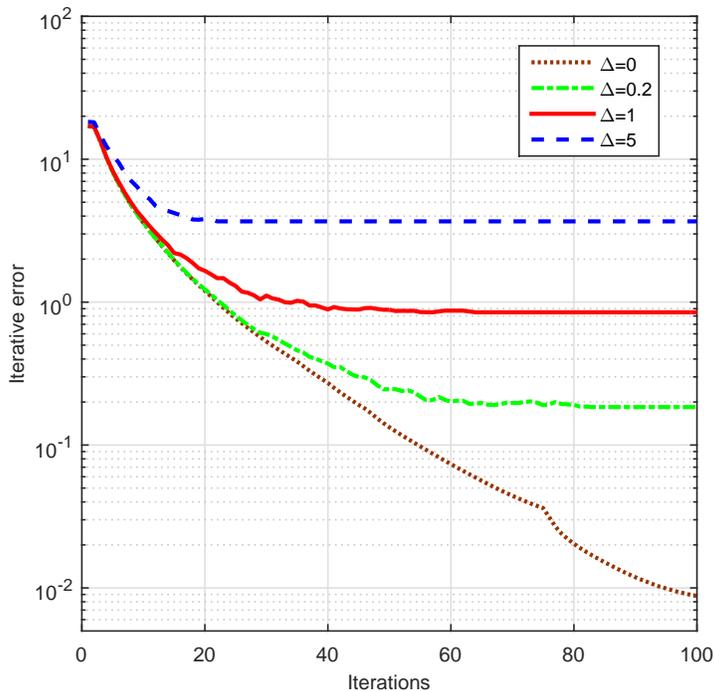}%
	\caption{The QC-ADMM for a distributed LASSO problem with $\rho=1$, $N=40$, $E=40$ and $\Delta\in\{0, 0.1, 0.5, 2.5, 10\}$.}
	\label{fig:QCADMMLASSOqr}
\end{figure}

{\it Quantization resolution:} Set $\rho=1$. Fig.~\ref{fig:QCADMMLASSOqr} is the simulation result of a network with $N=40$, $E=300$ and $\Delta\in\{0, 0.1, 0.5, 2.5 , 10\}$. Here $\Delta=0$ means that no quantization operation is placed on data communications, i.e., the C-ADMM is used to solve the LASSO problem. 

We observe that the consensus error becomes larger as $\Delta$ increases. This is not surprising as the higher the quantization resolution is, the more information is lost at each update, thus resulting in a higher consensus error. Meanwhile, the convergence time decreases when $\Delta$ increases, which can be seen from the upper bound on the number of iterations that guarantees the convergence of the QC-ADMM; that is, $\Omega$ in (\ref{eqn:OmegaUB}) decreases as $\Delta$ becomes larger. On the other hand, a larger $\Delta$ indicates a sparser quantization lattice which makes it easier for the QC-ADMM reach a convergence point.

{\it Algorithm parameter:} From Theorems \ref{thm:convergenceQCADMM} and \ref{thm:generalcase}, the upper bound on the consensus error increases with the algorithm parameter $\rho$. However, characterizing the effect of $\rho$ on the convergence time is very hard:
 $\rho$ not only affects the linear convergence rate $\eta$, but also involves in the upper bound on the number of iterations that guarantees the convergence of the QC-ADMM. Moreover, as shown in \cite{Wei:ADMMinDCO}, even though one can pick $\rho$ that maximizes the upper bound of $\eta$, the practical performance is usually suboptimal. Therefore, we only use a numerical example to study the effect of $\rho$ on the convergence time of the QC-ADMM.

\begin{figure}[htbp]
	\centering
	\includegraphics[width=0.7\textwidth]{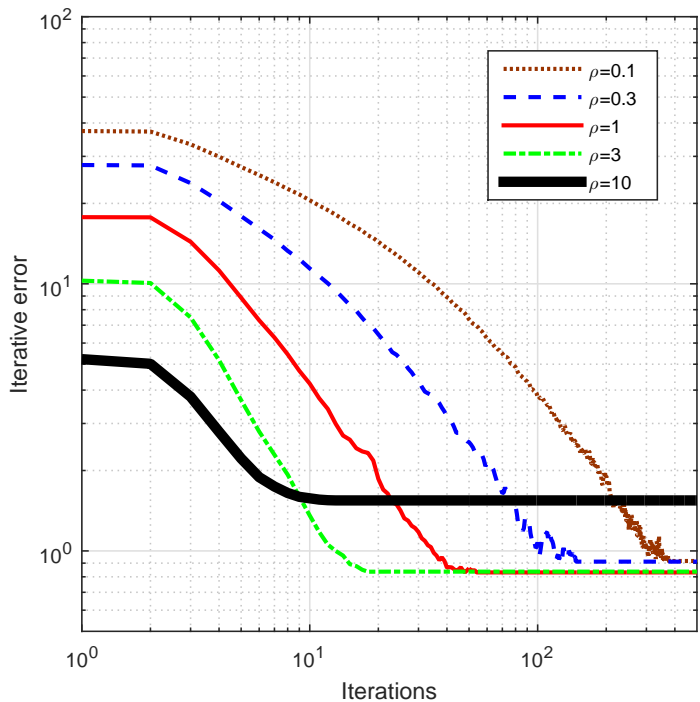}%
	\caption{The QC-ADMM for a distributed LASSO problem with $\Delta = 1$, $N=40$, $E=40$ and $\rho\in\{0.01, 0.1, 1,10\}$.}
	\label{fig:QCADMMLASSOrho}
\end{figure}

Set $\Delta=1$, $N=40$ and $E=300$. We apply the QC-ADMM to solving the above Lasso problem with $\rho \in\{0.01, 0.1, 1, 10\}$, and the result is presented in Fig.~\ref{fig:QCADMMLASSOrho}. In this example, the convergence time decreases as $\rho$ increases. Even though Theorem~\ref{thm:generalcase} indicates that a bigger $\rho$ results in a higher upper bound on the consensus error, the practical consensus error does not necessarily behave the same. That is, a bigger $\rho$ may lead a smaller consensus error in practice. 
%In this example, $m_g =0.0003 $ $M_g=72.14$.  $\tilde{\sigma}_{\min}=3.7386$ and $\sigma_{\max}(M_+)=8.06$. Therefore, we see. On the other hand. As an another observation, even though indicates that a higher $rho$ indicates a higher upper bound on consensus error, the practical consensus error does not behave similarly. That is, a higher $\rho$ may result in a smaller consensus error. 
%

\begin{figure}[htbp]
	\centering
	\includegraphics[width=0.7\textwidth]{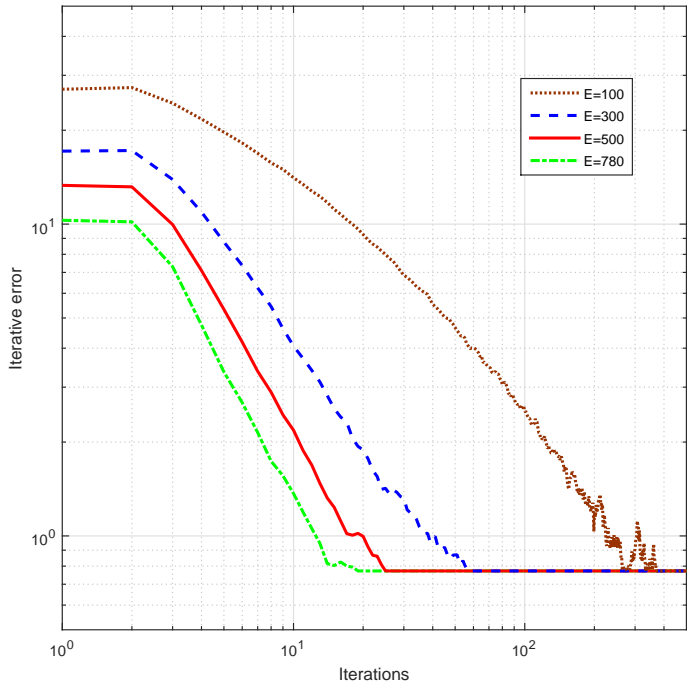}%
	\caption{The QC-ADMM for a distributed LASSO problem with $\Delta=1$, $\rho=1$, $N=40$, $E=\in\{100, 300, 500, 780\}$.}
	\label{fig:QCADMMLASSOgd}
\end{figure}

{\it Graph Density:} Fig.~\ref{fig:QCADMMLASSOgd} shows the result when $\Delta=1$, $\rho=1$, $N=40$, and $E\in\{100, 300, 500, 780\}$. The consensus error is the same for all $E$ while the convergence time decreases as $E$ becomes larger. Again, though Theorem \ref{thm:generalcase} indicates the upper bound for the consensus error increases in $E$ when $N$ is fixed, the practical consensus error need not necessarily perform the sam and can be much smaller. When $E$ increases with $N$ fixed, the average degree of the graph also increases. Then on the average, an agent can communicate with more agents at each update, thus resulting in a fast convergence. 
 
%\label{sec:simulations}
%\begin{align}
%\label{eqn:simulationfunction}
%\min_{\tilde{x}}\sum_{i=1}^N \frac{1}{2}\|U_i\tilde{x}-v_i\|_2^2+I_{\mathcal{X}_i}(\tilde{x})
%\end{align}

\section{Conclusion}
\label{sec:conclusion}
This paper proposes an efficient algorithm, the QC-ADMM, for multi-agent distributed optimization under the quantized communication constraint. We show that this algorithm can be derived from the standard ADMM by adding a quantization operation on $x^k$ immediately after the $x$-update together with proper initializations. While existing quantized ADMM approaches only apply to quadratic local objectives, the QC-ADMM can deal with more general objective functions, possibly non-smooth. Specifically, the QC-ADMM converges to a consensus within finite iterations under certain convexity conditions, which further enables us to derive a tight upper bound on the consensus error. Moreover, the proof idea provides a framework for convergence proof of a class of inexact updated algorithms.

Our approach also motivates future research directions: 

\begin{enumerate}
%\item We only consider the unbounded quantization scheme in this paper. There also exist applications that can be formulated in the form of problem (\ref{eqn:DistrOpt}) with bounded constraints on the variable, i.e., $\tilde{x}$ can only take values form some bounded set $\mathcal{X}$. Thus it is interesting to study the effect of the bounded quantization on our QC-ADMM.
\item We assume the quantized data communication between agents to be perfect in this paper. In practice, channel impairment may lead to imperfect transmissions. Moreover, the links between agents may fail and the topology of the network may vary. It is thus interesting to investigate how our algorithm performs in such settings. 
\item Recent work of \cite{Chang:ICADMM,Ling:DLADMM} proposes computationally efficient distributed ADMM algorithms that have linear convergence rates under certain conditions. We expect that the idea of this paper can also lead to quantized ADMM algorithms with significantly reduced computational complexity.
%\item linear constraint on the variables
%\item constrained ADMM
\end{enumerate}
%\section*{Acknowledgements}
%The authors would like to thank Professor Lixin Shen and Tiexing Wang for helpful discussions.
%\section*{Acknowledgments}
%The problem of this paper came into being when the first author visited the OSPAC group at the University of Minnesota and was later formulated when he returned to Syracuse University. The authors would like to thank Professor Zhi-Quan Luo and his group members for a happy visiting experience. The authors would also like to thank Professor Mingyi Hong for providing a related reference.
%$$ \text{cos}\theta = -\frac{6}{\sqrt{61}}$$
%$$ \text{sec}\theta = -\frac{\sqrt{61}}{6}$$
%$$ \text{cot}\theta = \frac{6}{5}$$

\bibliographystyle{IEEEbib}

\end{document}